\documentclass{article}
\usepackage[utf8]{inputenc}
\usepackage[utf8]{inputenc}
\usepackage{amsmath}
\usepackage{amsthm}
\usepackage{amssymb}
\usepackage{comment}
\usepackage{bbm}
\usepackage{comment}
\usepackage{float}
\usepackage{url}
\usepackage{geometry}
\usepackage{graphicx}
\usepackage[usenames,dvipsnames]{xcolor}
\usepackage{caption}
\usepackage{subcaption}
\usepackage{placeins}
\usepackage[colorinlistoftodos]{todonotes} 

\usepackage{tikz}

\usetikzlibrary{positioning}
\usetikzlibrary{arrows}
\usetikzlibrary{arrows.meta}
\usetikzlibrary{calc}

\definecolor{azzur}{RGB}{0,181,255}

\newcommand{\eps}{\varepsilon}

\title{Slow-fast dynamics in a neurotransmitter release model:\\ delayed response to a time-dependent input signal}
\author{Mattia Sensi$^{1,*}$, Mathieu Desroches$^1$ and Serafim Rodrigues$^{2,3}$\\[1em]
\small $^1$MathNeuro Team, Inria at Universit\'e C\^ote d'Azur, 06902 Sophia Antipolis cedex, France\\
\small $^2$MCEN research group, Basque Center for Applied Mathematics, 48009 Bilbao, Basque Country, Spain\\
\small $^3$Ikerbasque, The Basque Science Foundation, Bilbao, Basque Country, Spain.\\
\small $^*$Corresponding author. Email: \texttt{mattia.sensi@inria.fr}}
\date{}

\begin{document}

\maketitle

\begin{abstract}
    We propose a generalization of the neurotransmitter release model proposed in \emph{Rodrigues et al. (PNAS, 2016)}. We increase the complexity of the underlying slow-fast system by considering a degree-four polynomial as parametrization of the critical manifold. We focus on the possible transient and asymptotic dynamics, exploiting the so-called entry-exit function to describe slow parts of the dynamics. We provide extensive numerical simulations, complemented by numerical bifurcation analysis.
\end{abstract}


\section{Introduction}
Communication between neurons is at the core of brain activity and this is mainly organised by synapses. A majority of synapses are chemical, and involve the release of chemicals called \textit{neurotransmitters} in the small extracellular space called \textit{synaptic cleft} in between two neurons that are typically referred to as \textit{presynaptic neuron} and \textit{postsynaptic neuron}.

In~\cite{rodrigues2016time}, we succeeded in deriving a mathematical version of the biological model of neurotransmitter release (exocytosis) developed over many years by the 2013 Physiology or Medicine Nobel Prize laureate Thomas S{\"u}dhof and his collaborators, based (in particular) upon a fine understanding of the interactions between various protein complexes; see, e.g.,~\cite{sudhof2014,sudhof2009}. Our modeling approach was therefore to focus on key results obtained by S{\"u}dhof and collaborators on the dynamic interactions between two families of protein complexes, namely SNARE and SM, and how the dynamic changes of conformation between these two families of macromolecules affect the timing of release of the neurotransmitters. Namely, we started from the famous synaptic resources model by Tsodyks and Markram~\cite{markram1996,tsodyks1997}, which can only capture differential neurotransmitter release by hard-wiring a delay~\cite{wang2014} or adding stochastic terms to it~\cite{Senn2001}. This both leads to complex mathematical equations (delay-differential equations and stochastic differential equations, respectively) and does not really allow to model the biological processes underpinning the emergence of differential release of neurotransmitter.

Our approach in~\cite{rodrigues2016time} was different in at least two aspects. First, we wanted to put into equations the biological model put forth by S{\"u}dhof \emph{et al.}, namely the subcellular protein-protein interactions between SNARE and SM, responsible for the timing of neurotransmitter release. Second, we wanted to exploit the known presence of multiple timescales within these biochemical interactions in order to generate dynamical mechanisms giving rise to delays. Indeed, a known effect in the dynamics of systems displaying multiple timescales (usually referred to as \textit{slow-fast dynamical systems}), is to produce delayed transitions.\newline In terms of mechanistic behaviour of the model, the main add-on presented in~\cite{rodrigues2016time} is a two-dimensional slow-fast model where the variables $p_1$ and $p_2$ represent conformational states of the SNARE and SM protein complexes, respectively. This planar subsystem receives the electrical input formed by presynaptic spikes, and it transmits it with a controllable delay to the Tsodyks Markram synaptic resources model. This controllable delay is obtained as the result of a \textit{dynamic transcritical bifurcation}~\cite{krupa2001b} due to the timescale separation in the SNARE/SM component of the model and can be readily adjusted through parameter calibration and using an analytic formula referred to as \textit{entry-exit function}~\cite{neishtadt1987persistence,neishtadt1988persistence}; see already Section~\ref{sec:entryexit} for details.

The modeling approach proposed in~\cite{rodrigues2016time} was parsimonious and with an explicit mathematical structure based on slow-fast dynamical systems. Yet we could justify its derivation in biological terms by means of mass-action laws and inspiring from cell-cycle models~\cite{tyson1991}. Then, we used it to reproduce experimental data from both excitatory and inhibitory neurons. 

Our initial model from~\cite{rodrigues2016time} has a number of limiting aspects, in particular it only takes into account one part of the vesicle cycle. Therefore, it needs to be extended in order to capture the full cycle and better analyse synaptic transmission, for instance the factors that can make it malfunction. In the present work, we initialise an extension of our previous work on exocytosis on the theoretical ground by making a first attempt to include the endocycotic part of the cycle, where there is a re-uptake of the empty vesicle by the cellular membrane and then a refill of neurotransmitter molecules. Recent compelling experimental results~\cite{watanabe2013,watanabe2015} on this aspect of the vesicle cycle gives evidence of multiple pathways for the completion of the cycle after the exocytosis, with multiple timescales. This suggests a similar modeling strategy as the one we have used for the exocytotic part of the cycle.

To pave the way towards modeling the endocytotic as well as the exocytotic part of the vesicle cycle, our first attempt requires to increase the complexity of the model proposed in \cite{rodrigues2016time}, to be able to generate more complex behaviours such as postsynaptic responses with different amplitudes, akin to further differential release modes. To achieve this, we modify the fast nullcline, typically referred to as \textit{critical manifold}, changing the parabola to an M--shaped curve described as the graph of a fully factored fourth-degree polynomial.

We use of the entry-exit function to describe the passage of orbits close to a branch of the critical manifold, especially when this critical manifold loses normal hyperbolicity at a transcritical bifurcation point of the \textit{fast subsystem} (obtained by freezing the dynamics of the slow variable and considering it as a parameter), which is pivotal for the expected delayed transmission of the input signal. This enables us to 
first reproduce new types of outputs of the extended model, with both small- and large-amplitude oscillations of the postsynaptic membrane potential before going back to rest, upon time-dependent stimulus. These outputs correspond to additional patterns of differential neurotransmitter release, which qualitatively describe a first attempt to model multiple-timescale effects of the endocytotic part of the vesicle cycle. 

Furthermore, in absence of any time-dependent input and upon variation of one parameter of the extended model (which we denote $\alpha$), we analyse qualitatively and numerically the emergence of families of slow-fast limit cycles with intricate dynamical behaviour, due to the geometry of the underlying critical manifold. In particular, transverse intersections between a slow nullcline and a folded critical manifold is known to create \textit{canard cycles}~\cite{krupa2001} and we shall see that there are multiple families of such cycles in the extended system that we introduce here. Moreover, we analyse the convergence of the large relaxation cycles part of these families towards a heteroclinic orbit from the linear component of the critical manifold back to itself. Additionally, we compute the bifurcation structure of the model with respect to $\alpha$ by means of numerical continuation and use the software package \textsc{auto}~\cite{auto07p} to make such computations.

Amongst the rather large literature on multiple-timescale systems in neurodynamics, there is still little done on the effect of time-dependent input to such a slow-fast neural system. In particular, we already emphasized in~\cite{rodrigues2016time} the importance of multiple timescales on the generation of (possibly abnormal) delays in biological rhythms, aside the reliance upon time delays and stochastic noise in the model. In the present work, we promote a similar idea with a more complex geometry and describe the types of solutions obtained within this context.

The rest of the manuscript is organised as follows. In Section~\ref{sec:fullmodel}, we present the full model in detail and analyse its phase portrait, as well as, its slow-fast structure. The main technical tool underlying the delayed transmission of input displayed in the model is explained and analysed in Section~\ref{sec:entryexit}. Then, in Section~\ref{sec:numerics}, we perform a number of numerical exploration of the model in the two main dynamical regimes of interest: the stationary regime, where we simulate several scenarios of delayed response to input spikes; and the periodic regime, in absence of input stimulus, where we compute the full bifurcation structure of the model depending on the geometry of the quartic critical manifold. Finally, in Section~\ref{sec:conclusion}, we draw a few conclusions on this work and give some perspectives, both on further mathematical analysis of such model and on more detailed modelling that could inspire from it.

\section{The full neurotransmitter model and its slow-fast component}
\label{sec:fullmodel}

The full 6D system model derived and studied in \cite{rodrigues2016time} has the following form:
\begin{equation}\label{fullODE}
\begin{split}
    \Dot{p_1} &=  ( p_2-(a p_1+b) ) ( p_2-(\tilde{a} p_1+ \tilde{b}) )(\alpha -p_2)+V_{\text{in}}(t),\\ 
      \eps \Dot{p_2} &= p_2( p_1 -(\kappa_2 p_2^2 + \kappa_1 p_2 +\kappa_0) ),\\
      \Dot{d} &= (1-d)/\tau_D-dfp_2,\\
      \Dot{f} &= (f_0-f)/\tau_F+F(1-f)p_2,\\
      \Dot{g}_{\text{syn}} &= -g_{\text{syn}}/\tau_{\text{syn}}+\Bar{g}_{\text{syn}}dfp_2,\\
      C\Dot{v} &= -g_L(v-E_L)-g_{\text{syn}}(v-E_{\text{syn}}),
\end{split}
\end{equation}
where the overdot indicates the derivative with respect to the \textit{slow time $\tau$}, the variables $d$ and $f$ represent the minimal vesicle depletion model~\cite{markram1996,tsodyks1997}, $g_{\text{syn}}$ denotes the postsynaptic neuron conductance and $v$ the postsynaptic potential, which is the main experimentally observable quantity. As explained above, system~\eqref{fullODE} was introduced to describe synchronous, asynchronous and spontaneous release modes of neurotransmission between synaptic neurons; see~\cite{rodrigues2016time} for details.

The focus of the mathematical analysis in \cite{rodrigues2016time} is the study of the first two ODEs of system \eqref{fullODE}, representing a 2D system of ODEs in standard slow-fast form:
\begin{equation}\label{old}
\begin{split}
    p_1' &= \eps(p_2-(a p_1+b))(p_2-(\tilde{a} p_1+ \tilde{b}))(\alpha -p_2)+V_{\text{in}}(t)), \\ 
      p_2' &=p_2\left( p_1 -(\kappa_2 p_2^2 + \kappa_1 p_2 +\kappa_0) \right),
\end{split}
\end{equation}
where the $'$ indicates the derivative with respect to the fast time $t=\tau/\eps$, and $V_{\text{in}}(t)$ represents the stimulus received by the pre-synaptic neuron.

Both in \cite{rodrigues2016time} and here, we consider stimuli $V_{\text{in}}(t)$ of the form $V_{\text{in}}(t)=V \chi_{I}(t)$, with $V>0$ constant and $\chi_I(t)$ representing the characteristic function of the interval $I$, meaning
$$
\chi_I(t)=
\begin{cases}
1 & \text{ if } t\in I,\\
0 & \text{ otherwise.} 
\end{cases}
$$

To be able to capture a more dynamical scenario than in~\cite{rodrigues2016time}, we replace the parabolic nullcline of $p_2$ with a degree-four polynomial function, whose graph describes an M--shaped curve. The resulting system is the following slow-fast planar, with $p_1$ slow and $p_2$ fast:
\begin{equation}\label{new}
\begin{split}
    p_1' &= \eps(( p_2-(a p_1+b) ) ( p_2-(\tilde{a} p_1+ \tilde{b}) )(\alpha -p_2)+V_{\text{in}}(t)),\\ 
      p_2' &=p_2\left( p_1 -Q(c_1 p_2+r_1)(c_2 p_2+r_2)(c_3 p_2+r_4)(c_4 p_2+r_4) \right).
\end{split}
\end{equation}
We remark that we recover the original system \eqref{old} studied in \cite{rodrigues2016time} (up to renaming some parameters) if any two $c_i$ vanish. Moreover, as we showcase in our numerical simulations, when the dynamics is limited to the lower branch of the quartic, we recover a behaviour similar to something already observe in \cite{rodrigues2016time}. We refer to \cite{rodrigues2016time} for a detailed explanation of the parameters involved in system \eqref{old}, and we focus instead on the dynamics of system~\eqref{new}.

The \emph{critical manifold} of system \eqref{new} is defined from the slow-time version of the system, namely
\begin{equation}\label{new_slow_var}
\begin{split}
     \Dot{p_1}&=  ( p_2-(a p_1+b) ) ( p_2-(\tilde{a} p_1+ \tilde{b}) )(\alpha -p_2)+V_{\text{in}}(t),\\ 
     \eps \Dot{p_2} &=p_2\left( p_1 -Q(c_1 p_2+r_1)(c_2 p_2+r_2)(c_3 p_2+r_4)(c_4 p_2+r_4) \right),
\end{split}
\end{equation}
and taking the $\eps \rightarrow 0$ limit of the fast nullcline.

For ease of notation, we introduce the function
\begin{equation}\label{quartic}
    \Gamma:p_2\mapsto\;Q(c_1 p_2+r_1)(c_2 p_2+r_2)(c_3 p_2+r_4)(c_4 p_2+r_4).
\end{equation}
Then, the critical manifold $\mathcal{C}_0$ of \eqref{new} corresponds to the reunion of the $p_1$-axis with the graph of $\Gamma$, that is,
\begin{equation}\label{critmanif}
\mathcal{C}_0:=\{ p_2=0 \} \cup \{ p_1=\Gamma(p_2) \}.
\end{equation}
System \eqref{new} presents a transcritical point at the intersection of the two branches of the critical manifold, and we aim at exploiting this, combined with the slow passage close to $\{ p_2=0 \}$, to describe a delayed response to an input. This point is a transcritical bifurcation point of the \textit{fast subsystem} obtained as the $\eps=0$ limit of the fast-time system~\eqref{new}; see Fig.~\ref{phaseport} in which this point is labelled $TC$.

The system of interest for the present study is 2D, hence the only possible asymptotic behaviours are convergence towards a stable equilibrium or towards a stable limit cycle. However, we are also interested in deciphering the transient dynamics produced by system~\eqref{new_slow_var}, depending on the chosen values of the parameters as well as on the input stimulus $V_{\text{in}}$. In particular, we will investigate transients depending on parameters $\alpha$ and $r_1$.

The theoretical analysis of slow-fast systems like~\eqref{new_slow_var} is the realm of Geometric Singular Perturbation Theory (GSPT). Since its first conceptualization and the main results obtained by Fenichel~\cite{fenichel1979geometric}, GSPT has also proven to be remarkably useful to model real-world phenomena, which often naturally evolve on widely different timescales~\cite{bertram2017multi,hek2010geometric,jones1995geometric,wechselberger2020geometric}. In particular, it has been widely employed in the mathematical modelling of neural activity at various spatial and temporal scales; see e.g.~\cite{avitabile2017,desroches13,izhikevich00,moehlis06,rinzel86,Wechselberger13}.

\subsection{The quartic critical manifold}
For simplicity, and without loss of generality, we assume the following ordering on the zeros of $\Gamma$:
$$
0\leq -\dfrac{r_1}{c_1}<-\dfrac{r_2}{c_2}<-\dfrac{r_3}{c_3}<-\dfrac{r_4}{c_4}.
$$
The derivative of $\Gamma$ with respect to $p_2$ is a cubic polynomial function, characterized by
$$
\Gamma \bigg( -\dfrac{r_1}{c_1}\bigg) < 0, \quad \Gamma \bigg( -\dfrac{r_2}{c_2}\bigg) > 0, \quad \Gamma \bigg( -\dfrac{r_3}{c_3}\bigg) < 0, \quad \Gamma \bigg( -\dfrac{r_4}{c_4}\bigg) > 0.
$$
Hence, $\Gamma$ has exactly three critical points, respectively two negative minima and one positive maximum, located between consecutive zeros. The parameter $Q>0$ allows to modify the geometry of the critical manifold by controlling the location of the maximum and minima of $\Gamma$.

\begin{figure}[h!]\centering
\begin{tikzpicture}
    \node at (0,0){\includegraphics[width=0.49\textwidth]{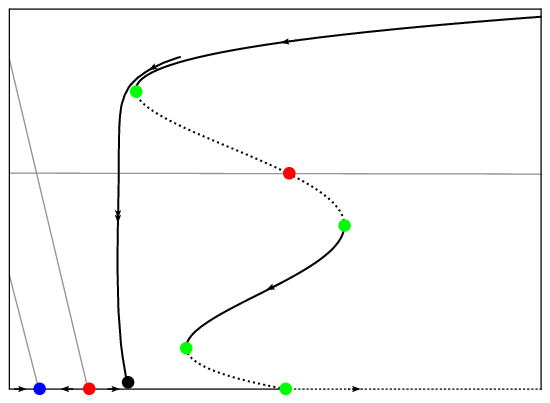}};	
    \node at (-3,-2.3) {\textcolor{blue}{$S$}};
    \node at (-2.25,-2.3) {\textcolor{red}{$U$}};
    \node at (0.4,0.75) {\textcolor{red}{$\tilde{U}$}};
    \node at (-1.75,-2.3) {\textcolor{black}{$L$}};
    \node at (0.25,-2.3) {\textcolor{green}{$TC$}};
    \node at (-0.6,-1.5) {\textcolor{green}{$SN_3$}};
    \node at (1.3,0) {\textcolor{green}{$SN_2$}};
    \node at (-1.2,1.5) {\textcolor{green}{$SN_1$}};
    \node at (2.25,2) {$p_2=\Gamma(p_1)$};
    \node at (3.6,-2.11) {$\rightarrow p_1$};
    \node at (3.48,-1.55) {$p_2$};
    \node at (3.25,-1.9) {$\uparrow$};
\end{tikzpicture}
\caption{Phase portrait of system \eqref{new}. The blue dot represents the stable equilibrium point $S=(-b/a,0)$, the red dots represent the unstable equilibria $U=(-\tilde{b}/\tilde{a},0)$ and $\tilde{U}=(\Gamma^{-1}(\alpha),\alpha)$, and the black dot represents the landing point $L$, located in a neighbourhood of $\{ p_2 = 0\}$, of an orbit (black curve) passing close to the upper branch of the quartic component of the critical manifold $\mathcal{C}_0$. Fast subsystem bifurcation points are represented by green dots: the transcritical point $TC=(\Gamma(0),0)$, which divides the stable (solid line) and unstable (dashed line) parts of $\{ p_2 = 0\}$; the saddle-node points $SN_i$, separating branches of the quartic $p_2 = \Gamma(p_1)$ with different stability. Grey lines correspond to the nullclines of $p_1$. The dynamics start in a neighbourhood of $S$, and the system is successfully excited if the input $V_{\text{in}}(t)$ is enough to ``jump'' to the right of $U$; otherwise, even after the input, the dynamics returns towards $S$ and stay in a vicinity of $\{ p_2 = 0\}$.}  \label{phaseport}
\end{figure}

Moreover, we assume $\Gamma(0)=Qr_1r_2r_3r_4\geq -\tilde{b}/\tilde{a}$, so that at $p_{1,TC}=\Gamma(0)$ the stability of the line $\{ p_2=0 \}$ changes from attracting to repelling, which corresponds to the transcritical bifurcation point of the fast subsystem; see Fig.~\ref{phaseport} for a visualization of the phase portrait of system~\eqref{new_slow_var}.

System~\eqref{new_slow_var} has one stable equilibrium and one unstable equilibrium located on the flat branch of the critical manifold; we denote them by $S$ and $U$, respectively. Depending on the relative positions of $S$, $U$ and $\Gamma(0)$, as well as on the initial conditions and the input stimulus $V_{\text{in}}(t)$, the system exhibits different transient behaviours, which will be explored in Section~\ref{sec:numerics}. 

The asymptotic behaviour can be summarized as follows. If $p_1(U)=-\tilde{b}/\tilde{a}$ is smaller than the landing off the higher branch of $\Gamma$, $L$, the dynamics approaches a stable limit cycle; if it is larger, after one or multiple inputs received in $p_2$, then the trajectory converges back to $S$.

\subsection{Equilibria and stability}

We refer again to Fig.~\ref{phaseport} for a visualization of the description carried out in this section. The critical manifold consists of two intersecting sets; recall from~\eqref{critmanif} that
$$\mathcal{C}_0=\{ p_2 = 0\} \cup \{ p_1 = \Gamma(p_2)\}.$$
The axis $\{ p_2=0 \}$ is a set of equilibria for the fast subsystem. On it, the only two equilibria of the full system are $S=(-b/a,0)$, stable, and $U=(-\tilde{b}/\tilde{a},0)$, unstable. These equilibria are given by the intersection of $\{ p_2=0 \}$ with the straight lines corresponding to the zeros of the first two factors of the $p_1$ nullcline.

Assuming that the lines $\{p_2=a p_1+b\}$ and $\{p_2=\tilde{a} p_1+\tilde{b}\}$ do not intersect the quartic curve $\{p_1=\Gamma(p_2)\}$, then the only remaining equilibrium is $\tilde{U}=(\Gamma^{-1}(\alpha),\alpha)$. The Jacobian of system \eqref{new} evaluated at $\tilde{U}$ has the following characteristic polynomial:
$$
\lambda^2-\alpha \Gamma'(\Gamma^{-1}(\alpha))\lambda +\eps\alpha (\alpha -(a\Gamma^{-1}(\alpha)+b))(\alpha-(\tilde{a}\Gamma^{-1}(\alpha)+\tilde{b}))=0.
$$
We are only interested in positive values of $\alpha$, as negative values would lead to unphysical solutions. Since, by assumptions, the two oblique components of the $p_1$-nullcline and $\Gamma$ do not intersect, the last term is strictly positive for any $\eps > 0$. Hence, as a consequence of Descartes' rule of signs, the eigenvalues will either have both positive or negative real parts. Namely, if $\Gamma'(\Gamma^{-1}(\alpha))>0$ (resp. $<0$), the real part of the eigenvalues will be positive (resp., negative), giving rise to another unstable (resp., stable) equilibrium.

If the equilibrium on $\Gamma$ is stable, orbits landing on the corresponding branch will be attracted to it. If it is unstable, orbits will not approach it, and will either converge back to a stable limit cycle -- a canard cycle if $\alpha$ is smaller than but close enough to the $p_2$ values at fold points of $\mathcal{C}_0$) -- or to $S$; the latter cases are illustrated on Figs.~\ref{cycle1phase} and~\ref{conveq1phase}, respectively.

The behaviour of an orbit near $\{p_2=0\}$ is key to the understanding of the dynamics of the system under scrutiny, namely the response of the system to (one or more) spike input(s). In particular, if we can control the input function $V_{\text{in}}(t)$, we know for which value of $p_1$ (or, equivalently, after how much time) $p_2$ will start to increase. To measure this, we make use of an important tool in GSPT, namely the \emph{entry-exit function}.

\section{Entry-exit function}
\label{sec:entryexit}
In this section, we focus on the application of the so-called entry-exit function \cite{kaklamanos2022entry,liu2000exchange,schecter2008exchange} (or delayed loss of stability \cite{neishtadt1987persistence,neishtadt1988persistence}) on system \eqref{new}, to describe the passage of orbits close to $\{ p_2=0 \}$. This function is remarkably useful to describe the permanence of a biological system near a manifold representing an absence of activity \cite{achterberg2022minimal,jardon2021geometric,jardon2021geometric2,rodrigues2016time}. For example, in epidemic models \cite{achterberg2022minimal,jardon2021geometric,jardon2021geometric2}, the entry-exit function can be exploited to compute the time between consecutive waves of an epidemic, characterized by a very small amount of infected individuals in the population.

The fast subsystem (or layer equation) of~\eqref{new} is the aforementioned system of ODEs after we let $\eps$ tend to $0$:
\begin{equation}\label{new_fast}
\begin{split}
    p_1' &= 0,\\ 
      p_2' &=p_2\left( p_1 -Q(c_1 p_2+r_1)(c_2 p_2+r_2)(c_3 p_2+r_4)(c_4 p_2+r_4) \right).
\end{split}
\end{equation}
System \eqref{new_fast} undergoes a transcritical bifurcation at the point $(\Gamma(0),0)$, labeled $TC$ in Figure \ref{phaseport}. Moreover, the axis $\{ p_2=0 \}$ changes its stability, from attracting to repelling, at $TC$, and the \textit{slow subsystem} on the same axis
$$
\Dot{p_1}= \alpha (a p_1+b)(\tilde{a} p_1+ \tilde{b}),
$$
ensures that $p_1$ is increasing to the right of the unstable equilibrium $U$, i.e. for $p_1 > -\tilde{b}/\tilde{a}$. Indeed, the slow subsystem (or reduced system) is defined as the $\eps=0$ limit of the slow-time system~\eqref{new_slow_var}, it is then a constrained system with the slow dynamics of the full system forced to evolve on the critical manifold. This means that on the $\{p_2=0\}$ component of $\mathcal{C}_0$, the slow system's dynamics, for variable $p_1$, is given by the above differential equation. Hence, the hypotheses necessary to apply the planar formulation of the entry-exit function are satisfied.

Assume the dynamics of system \eqref{new} starts close to the stable equilibrium $S$ at $(-b/a,\eps)$ and receives a kick strong enough to leap over the unstable equilibrium $U$ (if the kick is strong enough for the orbit to surpass the point $(\Gamma(0),0)$, we observe an immediate response).
The system then exhibits a classical planar entry-exit phenomenon \cite{achterberg2022minimal,de2008smoothness,de2016entry,jardon2021geometric}, as the $p_1$ axis changes stability at $\Gamma(0)$, being stable before and unstable afterwards, and $p_1$ increases on/near the branch of the critical manifold $\{ p_2=0 \}$. Namely, an orbit approaching the critical manifold between $U$ and $\Gamma(0)$ at a point with $p_1$ coordinate $p_{1,0}<\Gamma(0)$ will exit at a point $p_{1,1}>\Gamma(0)$ given by the solution of the following integral equation:
$$
\int_{p_{1,0}}^{p_{1,1}} \dfrac{x-\Gamma(0)}{(a x+b) (\tilde{a} x+ \tilde{b})} \text{d}x=0.
$$
This gives the implicit entry-exit relation as
\begin{equation}
    \label{eq:entry-exit}
a(\Gamma(0)\tilde{a}+\tilde{b})\log \bigg( \dfrac{\tilde{a}p_{1,0}+\tilde{b}}{\tilde{a}p_{1,1}+\tilde{b}} \bigg)=
\tilde{a}(\Gamma(0)a+b)\log \bigg( \dfrac{a p_{1,0}+b}{a p_{1,1}+b} \bigg).
\end{equation}
The higher branch of $\Gamma$ is ``accessible'' if and only if the exit point of an orbit entering the critical manifold close to $U$ is greater than the (local) maximum of $\Gamma$.
In this case, carefully balancing strength and duration of the kick, we can make sure that the first loop lands on the higher branch. Since the landing point from the both branches is to the right of $U$ (otherwise we would converge to the stable equilibrium $S$), the dynamics will then converge towards a limit cycle.

The piece of the limit cycle close to $\{ p_2 = 0 \}$ is completely characterized by the entry-exit relation: the limit cycle enters a neighbourhood of $\{ p_2 = 0 \}$ on the landing point from one of the two branches, and exits at a point given by relation \eqref{eq:entry-exit}.


\begin{figure}[h!]\centering
\begin{subfigure}[b]{0.49\linewidth}
 \centering
 \begin{tikzpicture}
   \node at (0,0){\includegraphics[width=\textwidth]{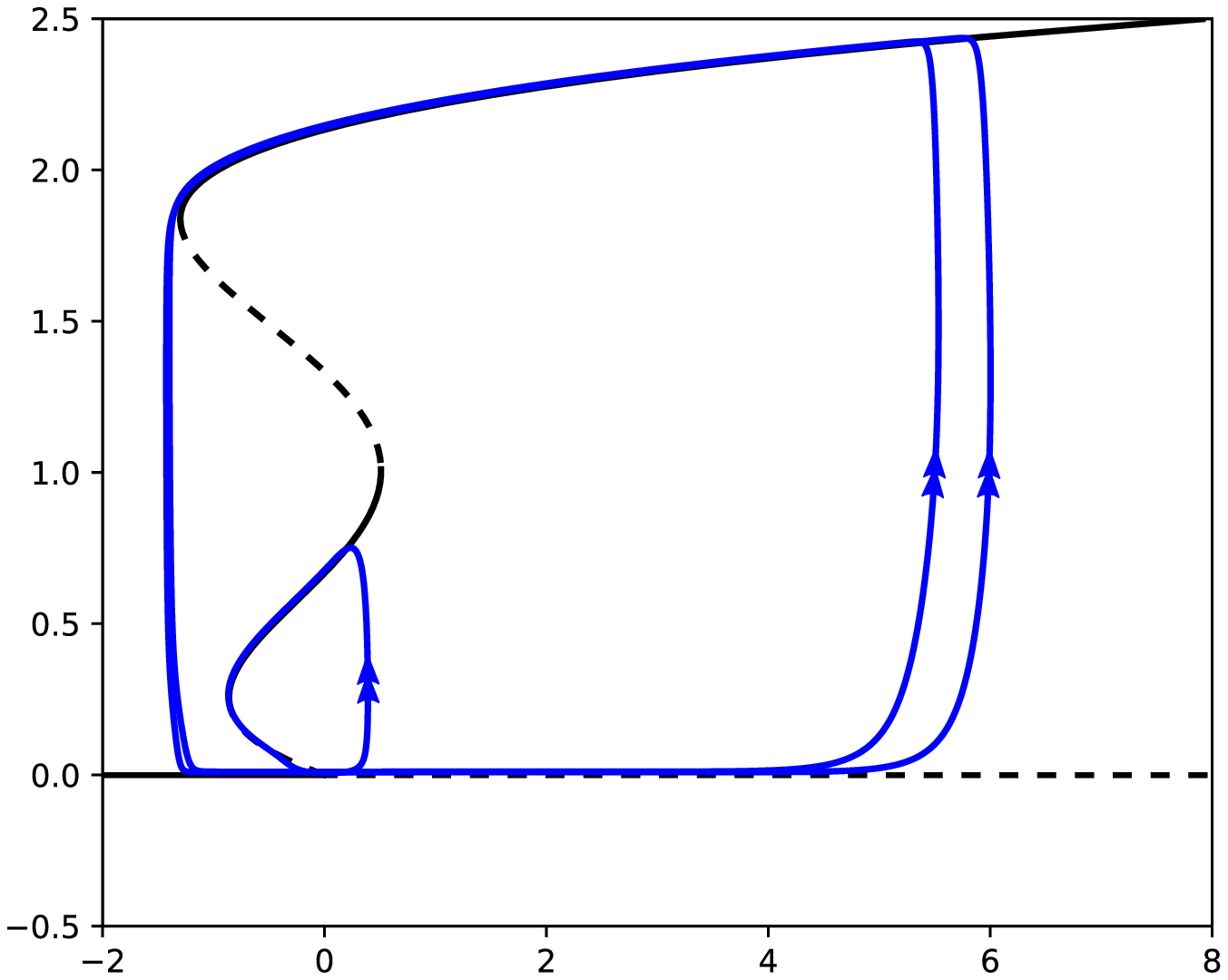}};	
   \node at (0,-3.25) {$p_1$};
   \node at (-3.75,0.25) {$p_2$};
	\end{tikzpicture} 
\caption{}\label{entrexsupper}
  \end{subfigure}
  \begin{subfigure}[b]{0.49\linewidth}
    \centering
  \begin{tikzpicture}
   \node at (0,0){\includegraphics[width=\textwidth]{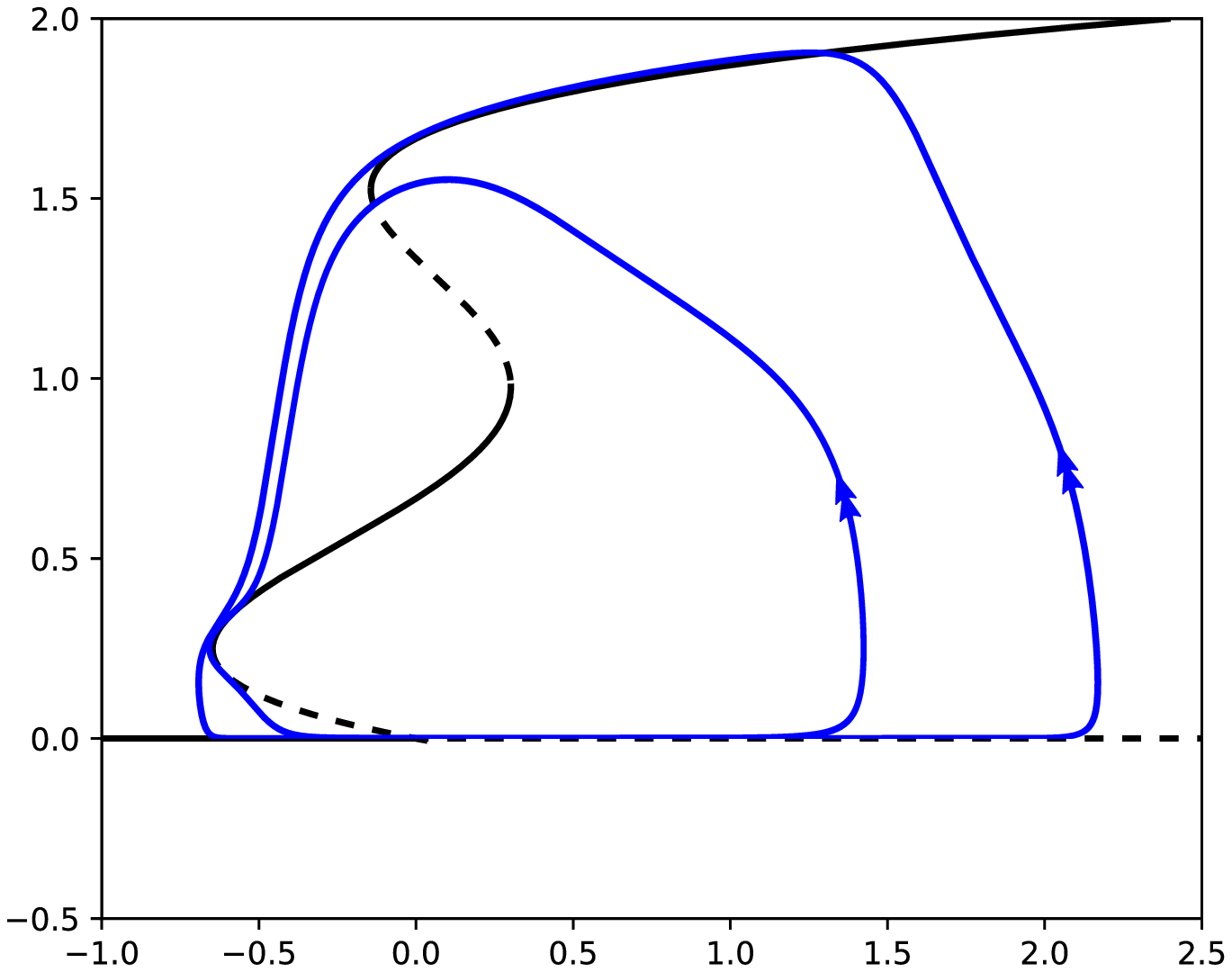}};	
   \node at (0,-3.25) {$p_1$};
   \node at (-3.75,0.25) {$p_2$};
   \end{tikzpicture}
\caption{}\label{entrexlower}
  \end{subfigure} 
\caption{Entry-exit phenomena for different values of $\alpha$, for both scenarios considered in Section~\ref{sec:numerics}. 
Panel (a) corresponds to the quartic \eqref{eqn:upper}; panel (b) to \eqref{eqn:lower}.
When $p_2<\alpha$, $p_1$ increases, so bigger values of $\alpha$ have bigger values for the entry-exit point in the entry-exit phenomena.
Notice how the entry points for the large orbits are very close, whereas the exit is greatly delayed.}
\label{entrexorbits}
\end{figure}

The parameter $\alpha$ has a crucial effect on the half-plane $\{p_2<\alpha\}$, which is the region where $p_1$ can grow, close to the quartic. In particular, this means that when a given trajectory falls off the quartic's leftmost fold point towards $\{p_2=0\}$, the landing point (which approximates the entry point in a neighbourhood of $\{p_2=0\}$) is increasing in $\alpha$. 

Hence, for bigger values of $\alpha$ we expect the entry to happen later, and consequently the exit to happen sooner, i.e. the trajectory will exit the neighbourhood of $\{p_2=0\}$ for smaller values of $p_1>\Gamma(0)$. After the exit, $p_1$ keeps increasing until the orbit intersects the line $\{p_2=\alpha\}$, then starts decreasing. Hence, the contribution of $\alpha$ to the maximum achieved on the limit cycle by $p_1$ is not monotone; see Figs.~\ref{entrexsupper} and~\ref{entrexlower} for a visualization of the role of $\alpha$ on the shape of the limit cycles.

The limit cycle given by the entry-exit from the ``falling point'' $L$ (recall Fig.~\ref{phaseport}) of the leftmost minimum of the quartic \emph{should} always be observed, when the equilibrium on the quartic is unstable, as the entry-exit formula is independent of $\alpha$.

It is possible to observe this limit cycle when the equilibrium on the quartic is stable too, as long as the cycle travels close to the stable branch of $\Gamma$ which does \emph{not} contain the equilibrium point. However, as discussed above, $\alpha$ has conflicting effects on the shape of the limit cycle. Smaller $\alpha$ means a smaller entry point, and consequently a greater exit point, but also less space for $p_1$ to grow to the maximum value it attains on the limit cycle.
\begin{figure}[h!]\centering
\begin{tikzpicture}
   \node at (0,0){\includegraphics[width=0.49\textwidth]{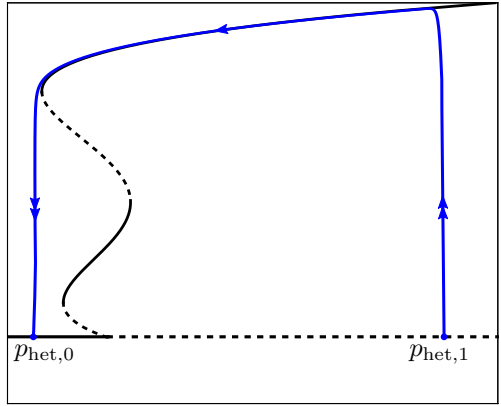}};	
   \node at (-2.5,-2) {$p_{\text{het},0}$};
   \node at (2.75,-2) {$p_{\text{het},1}$};
	\end{tikzpicture}
\caption{Limiting heteroclinic orbit as $\alpha$ tends to $0$. The $p_1$-axis consists of equilibria, attracting or repelling depending on the relative value of $p_1$ and $\Gamma(0)$. The ``landing'' point from the upper branch, $p_{\text{het},0}$, corresponding to the limit as $t \rightarrow +\infty$ of the limit heteroclinic orbit, gives the entry value to be used in \eqref{eq:entry-exit}; the other limit value of the heteroclinic orbit, as $t \rightarrow -\infty$, is given by the corresponding exit point, by a continuity argument as $\eps$ is small and $\alpha$ is close to $0$. The critical manifold $\mathcal{C}_0$ is represented with solid lines along its stable segments, and dashed lines along its unstable ones.}  \label{limalp}
\end{figure}

As $\alpha$ approaches $0$, orbits need to get increasingly closer to $\{p_2=0\}$ for $p_1$ to start increasing. In the $\alpha=0$ limit, the line $\{p_2=0\}$ consists of equilibria, which are attracting for $p_1<\Gamma(0)$ and repelling for $p_1>\Gamma(0)$. Hence there is no entry-exit behaviour for $\alpha=0$, and instead the limiting behaviour is a heteroclinic orbit between an unstable equilibrium of coordinates $(p_{\text{het},1})$ with $p_{\text{het},1}>\Gamma(0)$, and a stable equilibrium $(p_{\text{het},0})$ with $p_{\text{het},0}<\Gamma(0)$. Therefore, for $\alpha\neq0$ small, it takes longer and longer for orbits to pass close to $p_2=0$, resulting in limit cycles needing a increasingly high number of mesh intervals to be computed with orthogonal collocation as part of a numerical continuation procedure; we are using the software package \textsc{auto}~\cite{auto07p} for this. This explains the difficulties \textsc{auto} encounters in continuing the limit cycles close to $\alpha=0$. The limiting heteroclinic connection also includes a slow excursion close to the upper branch of the quartic curve $\{p_1=\Gamma(p_2)\}$. We sketch the shape of the limiting heteroclinic orbit as $\alpha \rightarrow 0$ in Figure \ref{limalp}.

In a sense, if we keep $\eps$ small but fixed, and we consider $\alpha$ as our perturbation parameter, the union of entry-exit strip on $\{ p_2 = 0 \}$ and the heteroclinic orbit to that same set, form a candidate orbit towards which the limit cycles tend, as $\alpha$ tends to $0$.
\FloatBarrier

\section{Numerical exploration of the full model}
\label{sec:numerics}
We will now study different scenarios of the delayed response to input in the quartic model~\eqref{new_slow_var}. These scenarios correspond to different configurations of the fold points of the quartic curve, namely the relative positions of the two left fold points. We will then consider two cases, depending on which one of the two is the leftmost one. We will also consider the transition case where both fold points are almost aligned. For each configuration, we will present numerical simulation of the full system and its delayed response to input spikes, as well as a detailed bifurcation structure of the limit cycles and equilibria (on $\Gamma$) of the $(p_1,p_2)$ system, with respect to $\alpha$. The values of the parameters in the last four ODEs of \eqref{fullODE} are selected among the ones used in \cite{rodrigues2016time}, namely: $\tau_D=200$,
$\tau_F=2500$,
$f_0=0.3$, 
$F=0.25$, 
$\tau_{\text{syn}}=20$,
$\Bar{g}_{\text{syn}}=0.4$, $C=0.196$, 
$g_L=1/220$,
$E_L=-55$ and
$E_{\text{syn}}=-57$.

To change the shape of the quartic, we only vary $r_1 \in [5,6.4]$, since this range is enough to substantially exchange the relative position of the outermost folds of the quartic, and consequently the bifurcation diagrams with respect to $\alpha$.

In all the figures, a single arrow represents the slow flow, whereas a double arrow represents the fast flow of the system. Moreover, blue represents stable limit cycles, and red unstable ones.

\subsection{Configuration 1: when the upper fold is smaller than the lower fold}
%
\begin{figure}[t!]\centering
\begin{subfigure}{.33\textwidth}
  \centering
\includegraphics[width=\textwidth]{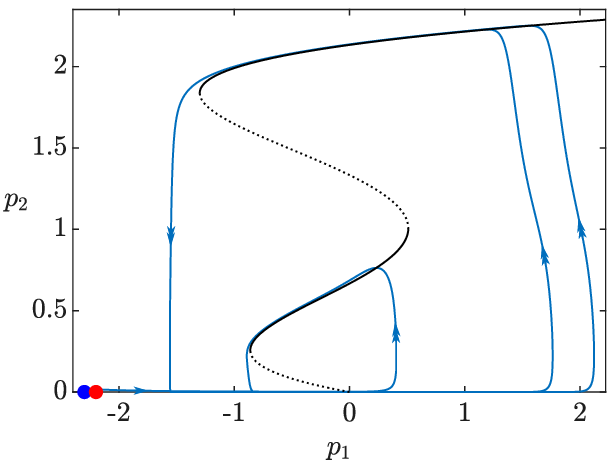}
\caption{}
  \label{cycle1phase}
\end{subfigure}%
\begin{subfigure}{.33\textwidth}
  \centering
  \includegraphics[width=\textwidth]{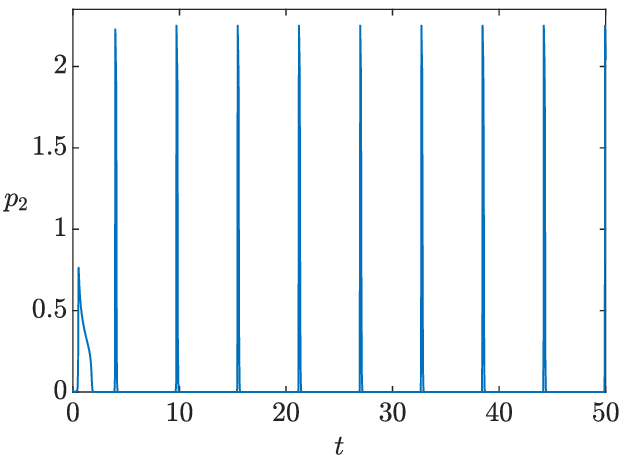}
\caption{}
  \label{cycle1p2}
\end{subfigure}
\begin{subfigure}{.33\textwidth}
  \centering
  \includegraphics[width=\textwidth]{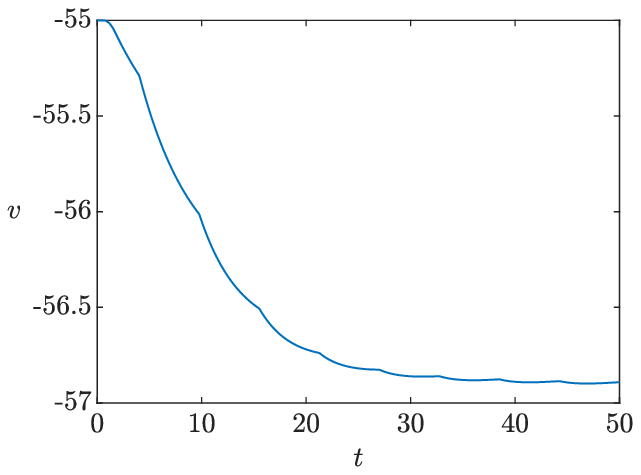}
\caption{}
  \label{cycle1v}
\end{subfigure}
\caption{(a): $(p_1,p_2)$ phase plane; (b): time series of $p_2$; (c): time series of $v$ in the full model. The blue and red dots represent $S$ and $U$, respectively. Parameter values are: $\eps=0.02$, $a=-1$, $b=-2.3$, $\tilde{a}=-1$, $\tilde{b}=-2.2$, $\alpha=0.22$,
$Q=0.05$, $c_1=-3$, $r_1=6.4$, $c_2=-3$, $r_2=4$, $c_3=-3$, $r_3=2$, $c_4=-3$, $r_4=0$. The starting point of the dynamics is slightly above the stable equilibrium $S$, namely $(p_1,p_2)=(-b/a,\eps)$. The input stimulus is the step function $V_{\text{in}}(t)=2700 \chi_{[0,0.04]}$ (where the indicator function is considered on $\mathbb{R}^+$). We observe two intermediate loops before the dynamics approaches the stable limit cycle. Notice that both equilibria on $\{p_2=0\}$ lie to the left of the falling point from the quartic $\Gamma$.}
\label{cycle1} 
\end{figure}
Here, we study the behaviour of system~\eqref{new} with the following choice for the quartic:
\begin{equation}\label{eqn:upper}
p_1=-(0.05(-3p_2+6.4)(-3p_2+4)(-3p_2+2)3p_2).
\end{equation}
The folds are located at $p_2 \approx 0.2565$ (corresponding to a local minimum), $p_2 \approx 1.00595$ (local maximum) and $p_2 \approx 1.8376$ (global minimum).
In order to observe (stable) limit cycles, both $S$ and $U$ need to be to the left of the ``falling point'' from the leftmost reaching branch of $\Gamma$, which depends on $\eps$ and is approximated by  $\min\Gamma(p_2)-\mathcal{O}(\eps^{2/3})$ (for the estimate as a trajectory jumps off a quadratic fold, we refer to~\cite{krupa2001}); see Fig.~\ref{cycle1} for such a configuration. 

\begin{figure}[t!]\centering
\begin{subfigure}{.33\textwidth}
  \centering
\includegraphics[width=\textwidth]{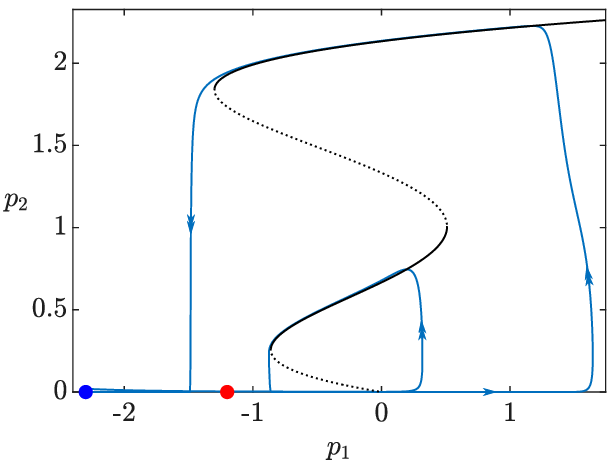}
\caption{}
  \label{conveq1phase}
\end{subfigure}%
\begin{subfigure}{.33\textwidth}
  \centering
  \includegraphics[width=\textwidth]{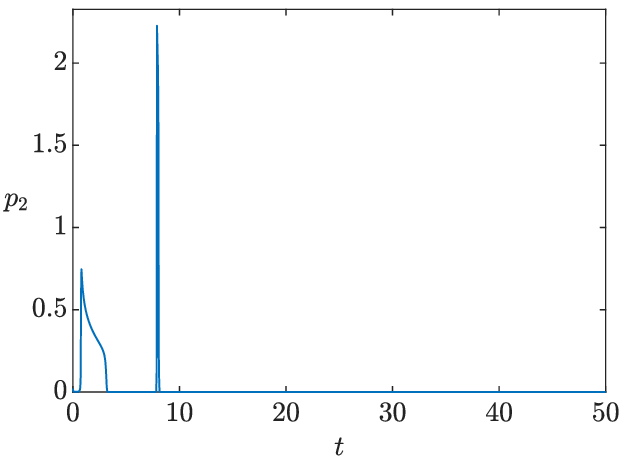}
\caption{}
  \label{conveq1p2}
\end{subfigure}
\begin{subfigure}{.33\textwidth}
  \centering
  \includegraphics[width=\textwidth]{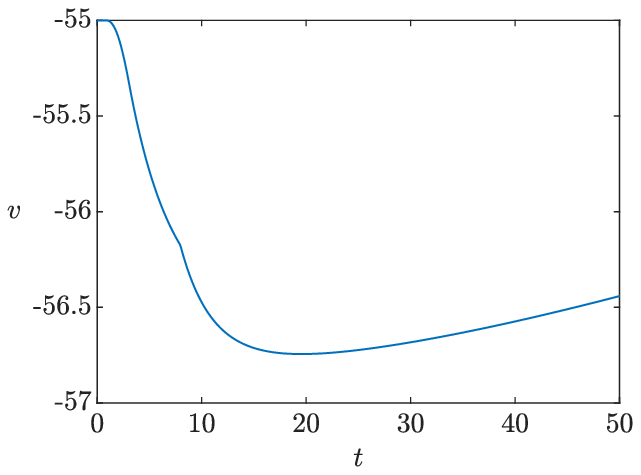}
\caption{}
  \label{conveq1v}
\end{subfigure}
\caption{(a): $(p_1,p_2)$ phase plane; (b): time series of $p_2$; (c): time series of $v$ in the full model. Parameter values are the same as in Figure \ref{cycle1}, except for $\tilde{b}=-1.2$, which places the unstable equilibrium $U$ between the two minima of $\Gamma$. The starting point of the dynamics is slightly above the stable equilibrium $S$, namely $(p_1,p_2)=(-b/a,\eps)$. The input stimulus is the step function $V_{\text{in}}(t)=2700 \chi_{[0,0.04]}$ (where the indicator function is considered on $\mathbb{R}^+$). We observe a small loop, followed by a larger one, before the dynamics converge back to the stable equilibrium $S$.}
\label{conveq1} 
\end{figure}

Alternatively, whether both $S$ and $U$ lie between the ``falling point'' $L$ and $\Gamma(0)$ or only $U$ does, such as in Fig.~\ref{conveq1}, as long as the upper branch of $\Gamma$ is accessible to the dynamics, after a finite amount of spikes in the time series of $p_2$, orbits converge back to the stable equilibrium $S$.

\begin{figure}[h!]\centering
\begin{subfigure}{.33\textwidth}
  \centering
\includegraphics[width=\textwidth]{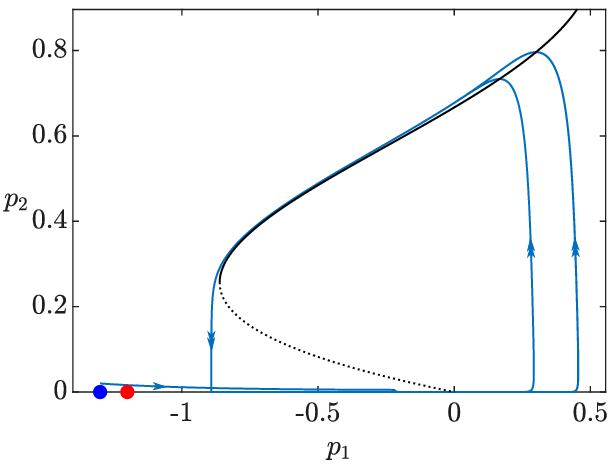}
\caption{}
  \label{smallcycle1phase}
\end{subfigure}%
\begin{subfigure}{.33\textwidth}
  \centering
  \includegraphics[width=\textwidth]{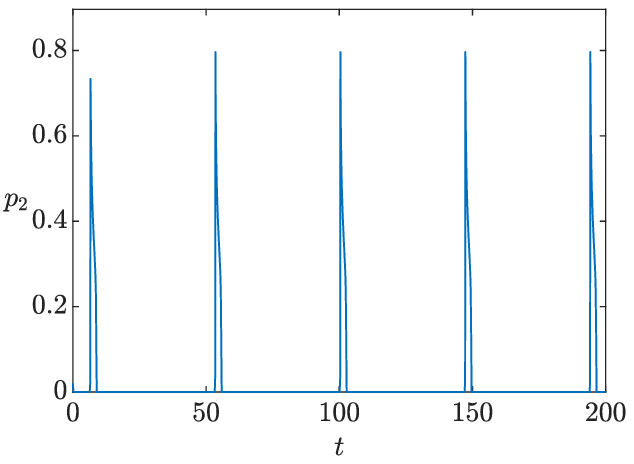}
\caption{}
  \label{smallcycle1p2}
\end{subfigure}
\begin{subfigure}{.33\textwidth}
  \centering
  \includegraphics[width=\textwidth]{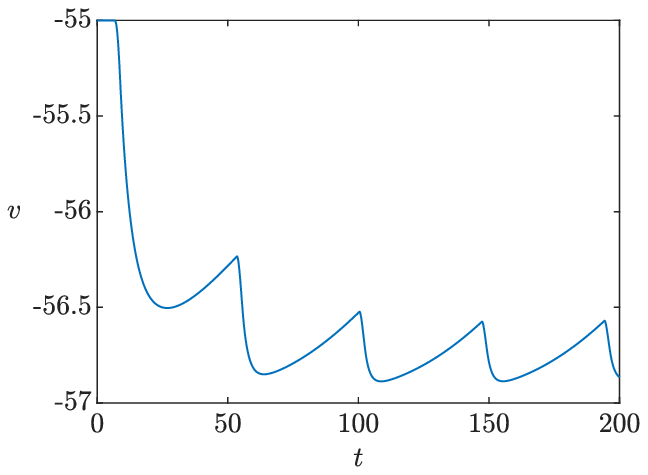}
\caption{}
  \label{smallcycle1v}
\end{subfigure}
\caption{(a): $(p_1,p_2)$ phase plane; (b): time series of $p_2$; (c): time series of $v$ in the full model. The blue and red dots represent $S$ and $U$, respectively. Parameter values are the same as in Figure \ref{cycle1}, except for $\alpha=0.05$, $b=-1.3$ and $\tilde{b}=-1.2$. The starting point of the dynamics is slightly above the stable equilibrium $S$, namely $(p_1,p_2)=(-b/a,\eps)$. The input stimulus is the step function $V_{\text{in}}(t)=1350 \chi_{[0,0.04]}$ (where the indicator function is considered on $\mathbb{R}^+$). We observe some intermediate loops on the lower branch before the dynamics approaches the stable limit cycle. Notice that both equilibria on $p_2=0$ lie to the right of the falling point from the quartic $\Gamma$ (compare with Figs.~\ref{cycle1} and~\ref{conveq1}).}
\label{smallcycle1} 
\end{figure}

Another possibility is to be trapped in a smaller loop on the lower branch, even if this condition is not satisfied; see Fig.~\ref{smallcycle1} for such a configuration.
\begin{figure}[h!]\centering
\begin{subfigure}{.33\textwidth}
  \centering
\includegraphics[width=\textwidth]{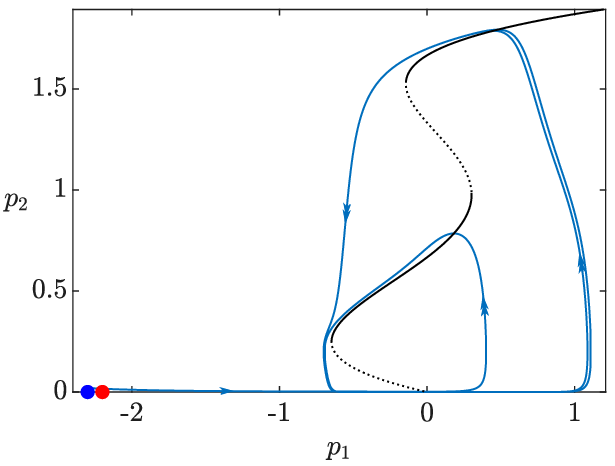}
\caption{}
  \label{cycle1SWphase}
\end{subfigure}%
\begin{subfigure}{.33\textwidth}
  \centering
  \includegraphics[width=\textwidth]{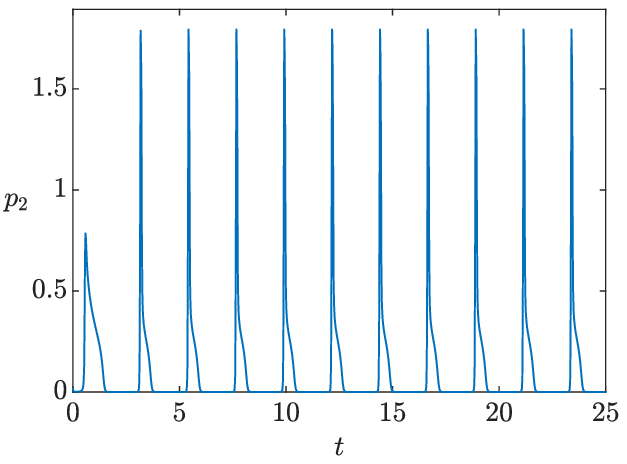}
\caption{}
  \label{cycle1SWp2}
\end{subfigure}
\begin{subfigure}{.33\textwidth}
  \centering
  \includegraphics[width=\textwidth]{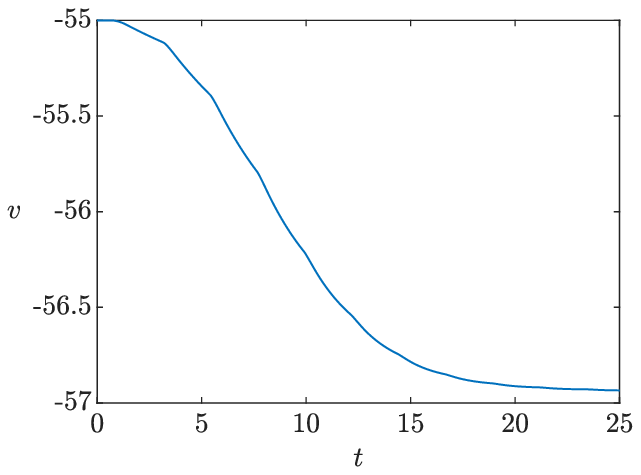}
\caption{}
  \label{cycle1SWv}
\end{subfigure}
\caption{(a): $(p_1,p_2)$ phase plane. (b): time series of $p_2(t)$. (c): time series of $v$ in the full system~\eqref{fullODE}. The blue and red dots represent $S$ and $U$, respectively. Values of the parameters: $\eps=0.02$, $a=-1$, $b=-2.3$, $\tilde{a}=-1$, $\tilde{b}=-2.2$, $\alpha=0.22$,
$Q=0.05$, $c_1=-3$, $r_1=5$, $c_2=-3$, $r_2=4$, $c_3=-3$, $r_3=2$, $c_4=-3$, $r_4=0$. The starting point of the dynamics is slightly above the stable equilibrium $S$, namely $(p_1,p_2)=(-b/a,\eps)$. The ``kick'' is the step function $V_{\text{in}}(t)=2700 \chi_{[0,0.04]}$. We observe one intermediate loop before the dynamics approaches the stable limit cycle, which exhibits visible slow passages on both branches of the quartic $\Gamma$. Notice that both equilibria on $p_2=0$ lie to the left of the quartic $\Gamma$.}
\label{cycle2} 
\end{figure}
%

\subsubsection{Bifurcation structure in $\alpha$}
We are interested in the behaviour of system~\eqref{new} as we vary the value of $\alpha$. In particular, we remark that in a $(\alpha,p_1)$ bifurcation diagram, the curve of equilibria for $p_1>0$ will be represented as $p_1=\Gamma(\alpha)$.

\begin{figure}[h!] \centering
  \begin{subfigure}[b]{0.49\linewidth}
    \centering
    \begin{tikzpicture}[thick,scale=0.8, every node/.style={scale=0.8}]
   \node at (0,0){\includegraphics[width=\linewidth]{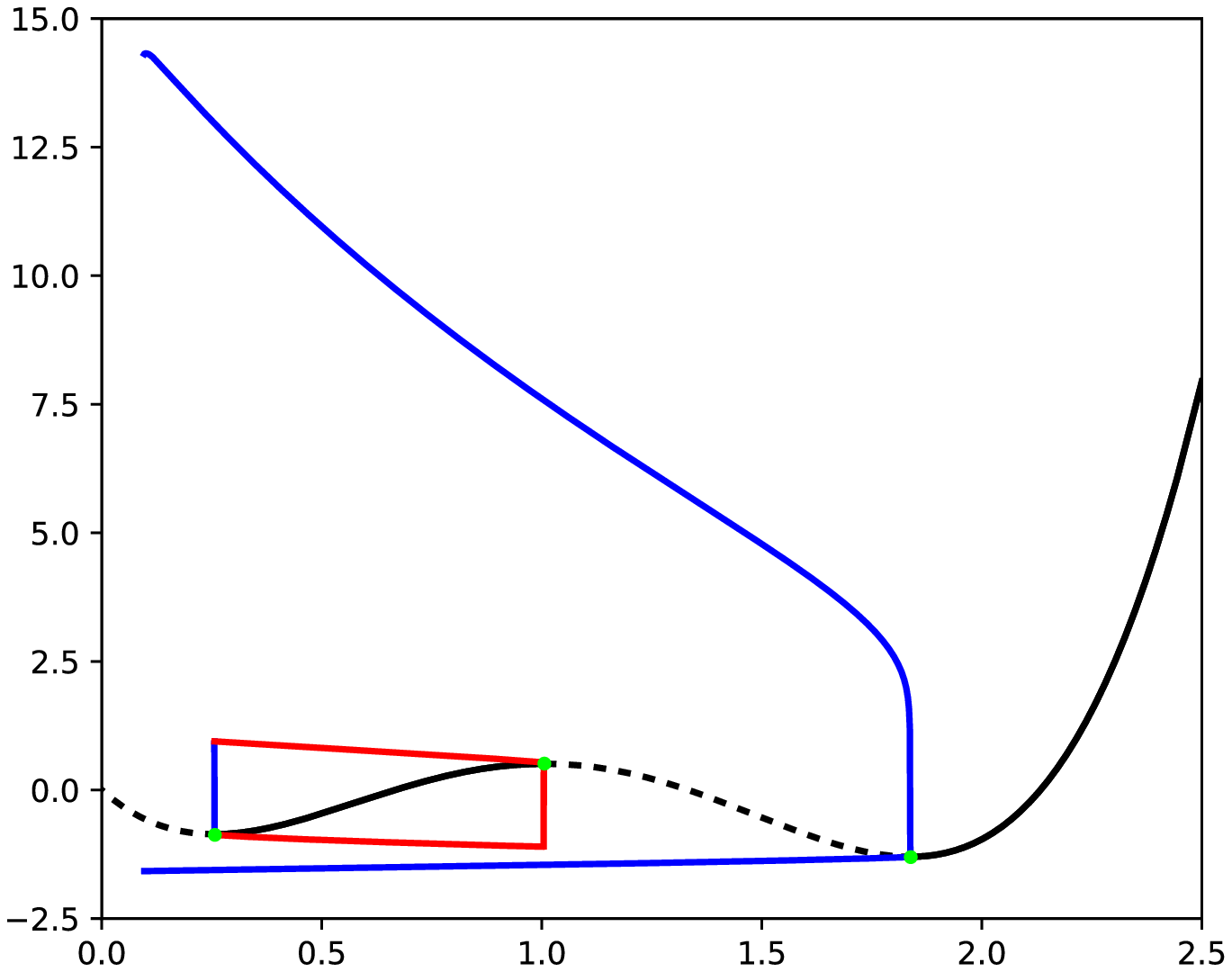}};	
  \node at (0,-3.25) {$\alpha$};
   \node at (-3.5,0) {$p_1$};
    \node at (2,-2.5) {\textcolor{green}{$H_1$}};
   \node at (-0.5,-1.5) {\textcolor{green}{$H_2$}};
    \node at (-2.6,-1.9) {\textcolor{green}{$H_3$}};
	\end{tikzpicture} 
    \caption{} 
    \label{bifupper:total} 
  \end{subfigure}
  \begin{subfigure}[b]{0.49\linewidth}
    \centering
    \begin{tikzpicture}[thick,scale=0.8, every node/.style={scale=0.8}]
   \node at (0,0){\includegraphics[width=\linewidth]{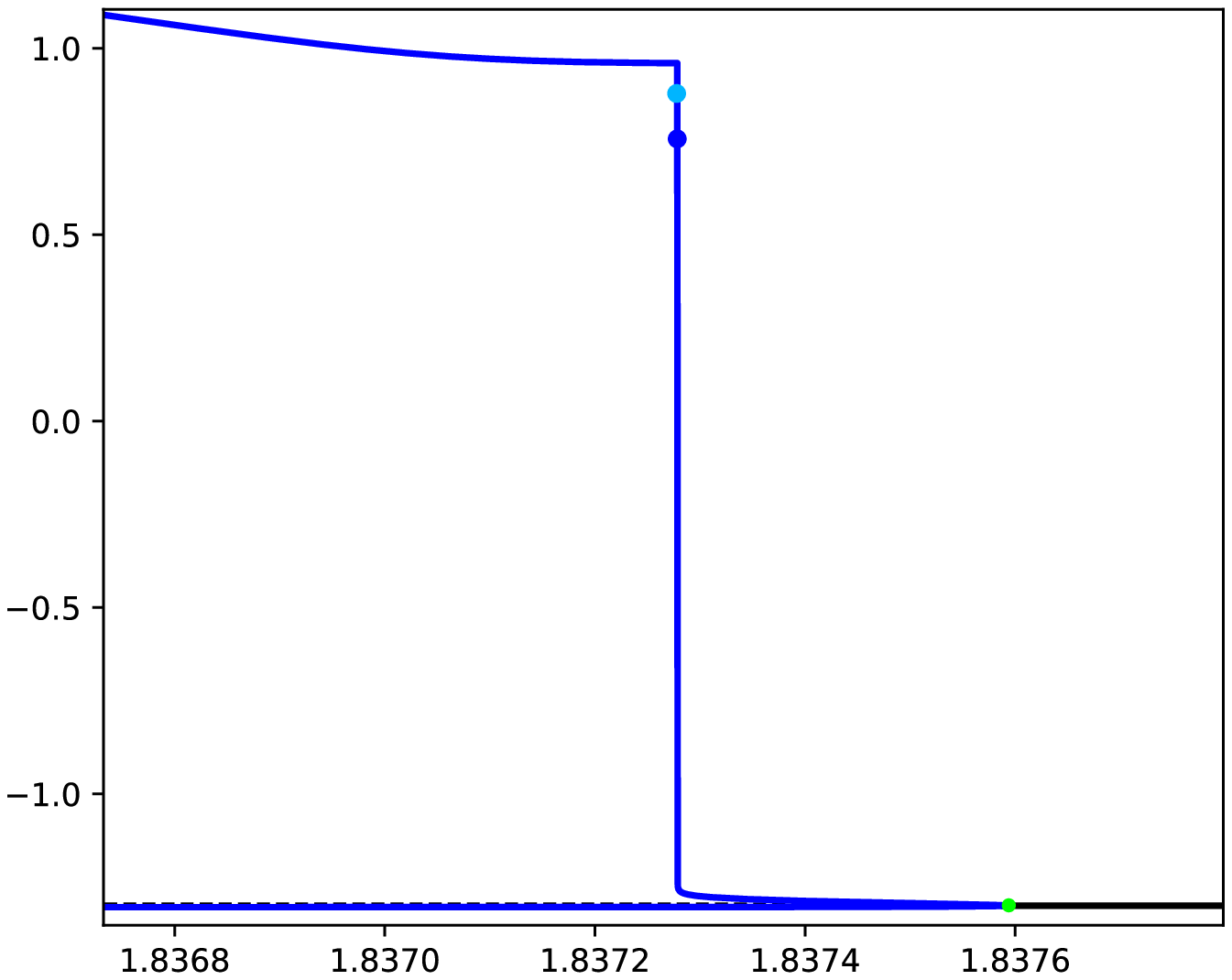}};	
    \node at (0,-3.25) {$\alpha$};
   \node at (-3.5,0) {$p_1$};
    \node at (2.2,-2.3) {\textcolor{green}{$H_1$}};
   \node at (0.5,2.2) {\textcolor{azzur}{$2$}};
  \node at (0.5,1.8) {\textcolor{blue}{$1$}};
	\end{tikzpicture} 
    \caption{} 
    \label{bifupper:hopf1}
  \end{subfigure} 
  \begin{subfigure}[b]{0.49\linewidth}
    \centering
    \begin{tikzpicture}[thick,scale=0.8, every node/.style={scale=0.8}]
   \node at (0,0){\includegraphics[width=\linewidth]{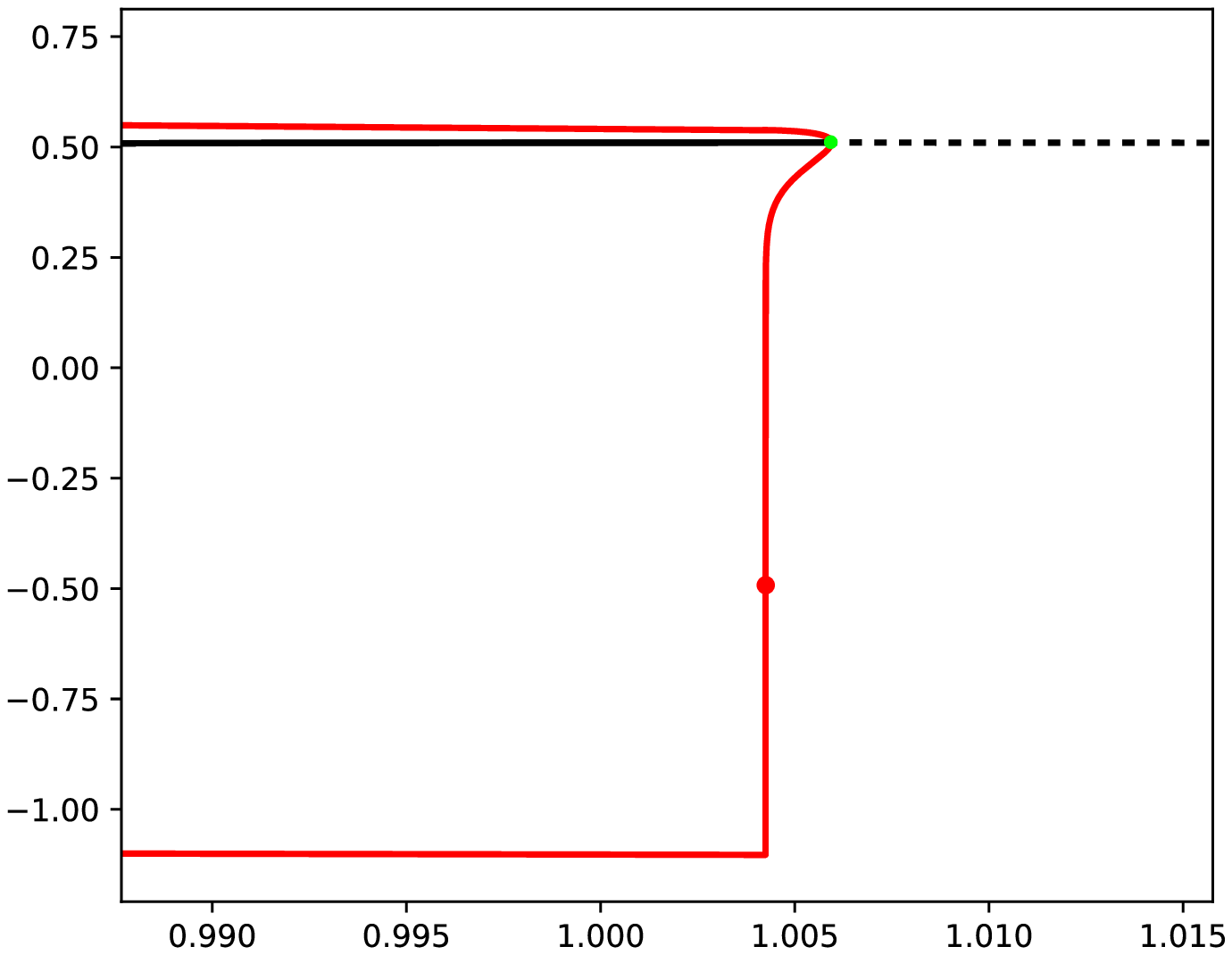}};	
   \node at (0,-3.25) {$\alpha$};
   \node at (-3.5,0.2) {$p_1$};
    \node at (1.4,2.2) {\textcolor{green}{$H_2$}};
   \node at (1.05,-0.75) {\textcolor{red}{$3$}};
	\end{tikzpicture} 
    \caption{} 
   \label{bifupper:hopf2} 
  \end{subfigure}
  \begin{subfigure}[b]{0.49\linewidth}
    \centering
    \begin{tikzpicture}[thick,scale=0.8, every node/.style={scale=0.8}]
   \node at (0,0){\includegraphics[width=\linewidth]{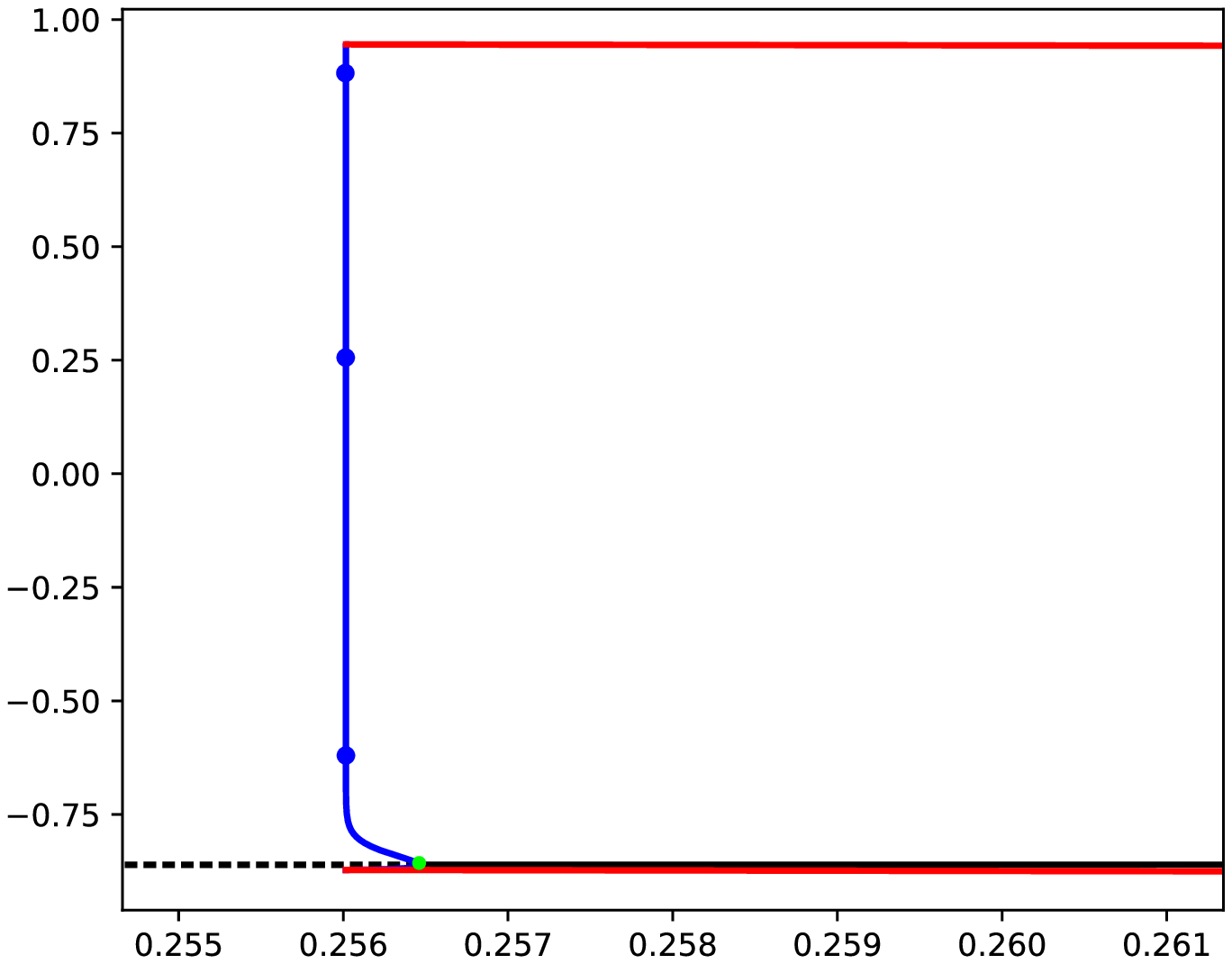}};	
   \node at (0,-3.25) {$\alpha$};
   \node at (-3.5,0.2) {$p_1$};
    \node at (-1.25,-2.1) {\textcolor{green}{$H_3$}};
   \node at (-1.85,-1.7) {\textcolor{blue}{$6$}};
   \node at (-1.85,0.65) {\textcolor{blue}{$5$}};
   \node at (-1.85,2.3) {\textcolor{blue}{$4$}};
	\end{tikzpicture} 
    \caption{} 
   \label{bifupper:hopf3} 
 \end{subfigure} 
  \caption{(a): Bifurcation diagram of the $(p_1,p_2)$ system with respect to $\alpha$, showing branches of equilibria represented by their $p_1$ value, and branches of limit cycles for which both maximum and minimum of $p_1$ along the cycles, are plotted. Solid black lines correspond to stable equilibria, while dashed black lines correspond to unstable ones. Stable limit cycle branches are shown in blue, and unstable ones in red. (b), (c), (d): zoomed view of the canard-explosive branches emanating from each of the three Hopf bifurcation shown (black dots) in panel (a), from right to left. The numbered dots correspond to canard orbits shown in Fig.~\ref{canardsupper}.}
  \label{bifupper1}
\end{figure}

For $\alpha$ greater than the first fold of the quartic (the rightmost in Fig.~\ref{bifupper1}), the corresponding equilibrium on the quartic is stable. At $\alpha=\alpha(H_1)$, which indicates the greatest value of $\alpha$ at which a Hopf bifurcation occurs, we observe a \textit{canard explosion}~\cite{krupa2001}; this branch extends to values of $\alpha$ close to $0$. 
However, in this region $\alpha$ qualitatively acts as an additional perturbation parameter in the part of the cycle close to $p_2=0$, and this introduces numerical issues, as the period of limit cycles grows unboundedly when they converge to a heteroclinic orbit. This leads to non-convergence of the continuation algorithm in \textsc{auto} when approaching $\alpha=0$, although we would analytically expect the continuation branch to reach $\alpha=0$ at a finite value of $p_1$.

The other two Hopf points on the quartic are connected by a branch of limit cycles, exhibiting canard explosions at both extrema. As we will see in the following bifurcation analysis, changing the value of $r_1$ will alter the connection between the Hopf points, from $H_2$-$H_3$ to $H_1$-$H_2$.

\begin{figure}[h!]\centering
\begin{subfigure}[b]{0.49\linewidth}
    \centering
\begin{tikzpicture}[thick,scale=0.8, every node/.style={scale=0.8}]
   \node at (0,0){\includegraphics[width=\linewidth]{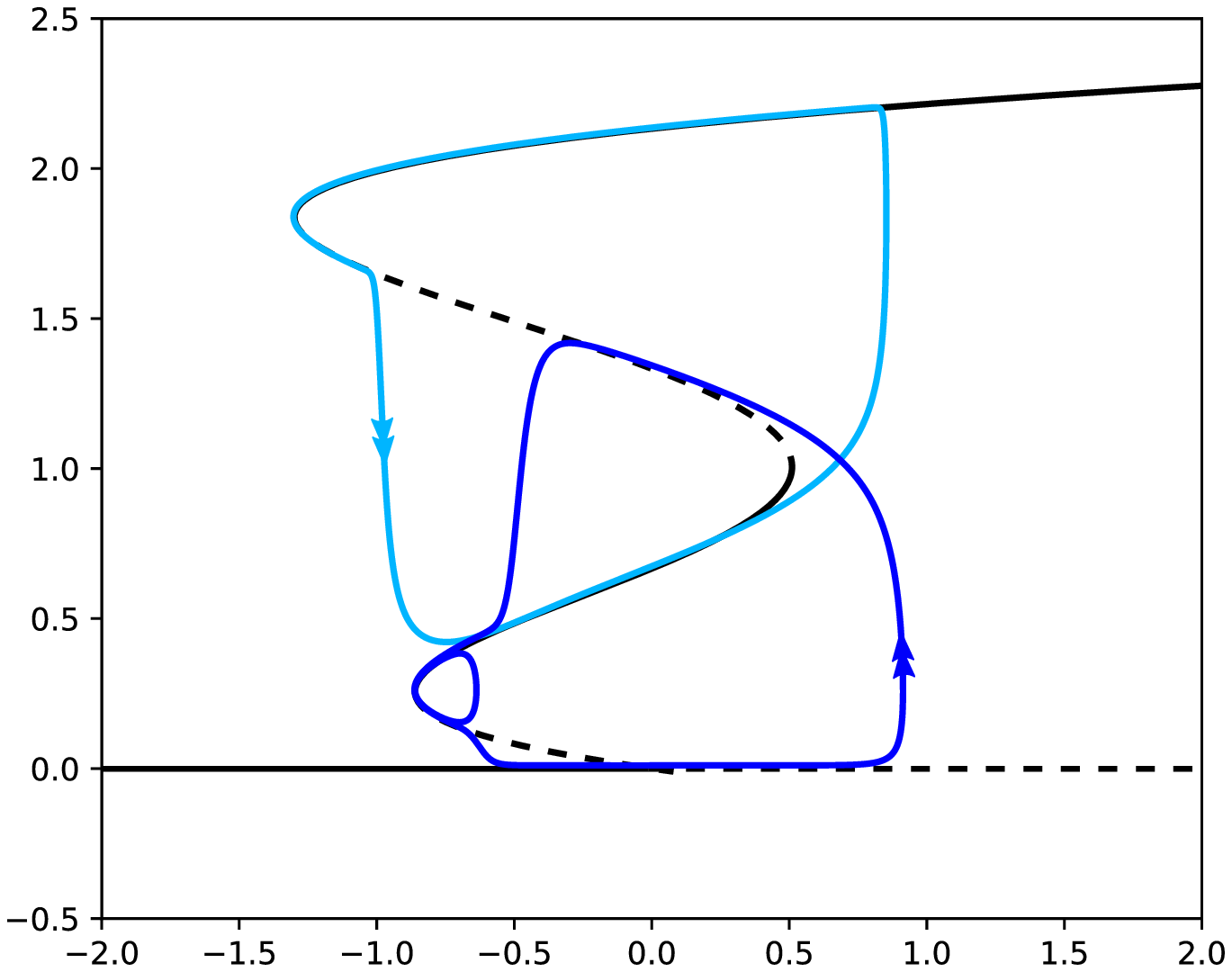}};	
   \node at (0,-3.25) {$p_1$};
   \node at (-3.5,0.25) {$p_2$};
   \node at (-1.3,2) {\textcolor{azzur}{$2$}};
   \node at (-0.6,-1.25) {\textcolor{blue}{$6$}};
   \node at (2,-1.25) {\textcolor{blue}{$4$}};
	\end{tikzpicture}  
    \caption{} 
    \label{canupper1} 
  \end{subfigure}
  \begin{subfigure}[b]{0.49\linewidth}
    \centering
\begin{tikzpicture}[thick,scale=0.8, every node/.style={scale=0.8}]
   \node at (0,0){\includegraphics[width=\linewidth]{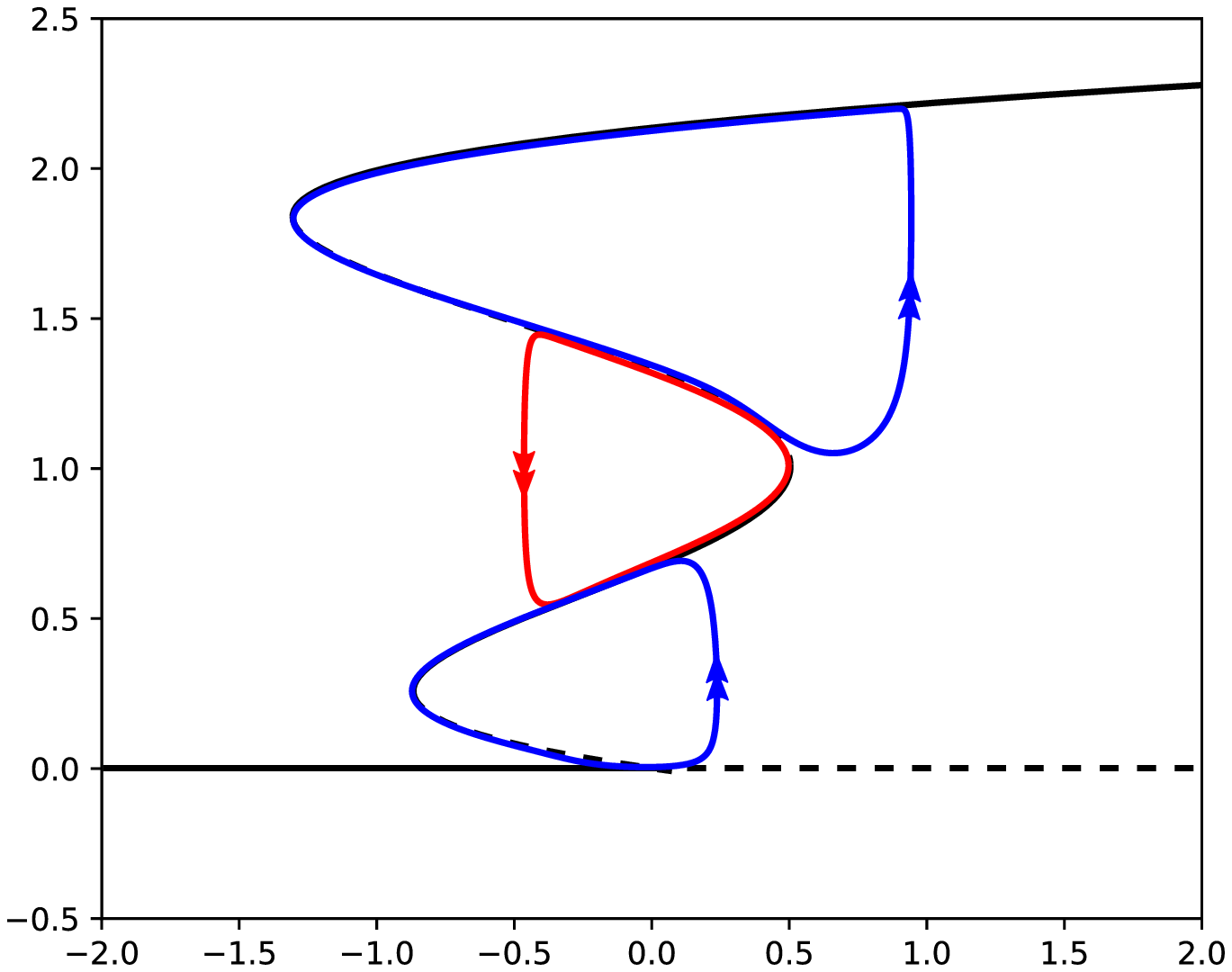}};	
   \node at (0,-3.25) {$p_1$};
   \node at (-3.5,0.25) {$p_2$};
  \node at (2,1.25) {\textcolor{blue}{$1$}};
   \node at (-0.8,0) {\textcolor{red}{$3$}};
   \node at (0.8,-1.25) {\textcolor{blue}{$5$}};
	\end{tikzpicture}  
    \caption{} 
    \label{canupper2} 
  \end{subfigure} 
\caption{Selection of canard cycles that exist along the branches of limit cycles shown in Fig.~\ref{bifupper1} (b), (c) and (d).}  
\label{canardsupper} 
\end{figure}

In Fig.~\ref{bifupper1} panels (b-d), we highlighted with dots, and numbered from 1 to 6, points on the bifurcation diagram corresponding to a selection of canard orbits. These orbits are then plotted in the phase plane in Fig~\ref{canardsupper}. The cycles are numbered in decreasing values of $\alpha$, although clearly the cycles belonging to the same canard explosive branch correspond to values of $\alpha$ extremely close.
%
%

 
\subsection{Configuration 2: when the lower fold smaller than the upper fold}
In this section, we numerically explore the behaviour of system \eqref{new} with the following choice for the quartic:
\begin{equation}\label{eqn:lower}
p_1=-(0.05(-3p_2+5)(-3p_2+4)(-3p_2+2)3p_2).
\end{equation}
The folds are at $p_2 \approx 0.24875$ (corresponding to the global minimum), $p_2 \approx 0.976512$ (local maximum) and $p_2 \approx 1.52474$ (local minimum).

If we compare Fig.~\ref{cycle2} with Fig.~\ref{cycle1}, we notice that in this case the passage close to the lower branch of the quartic gives rise to another slow section in each spike.

\subsubsection{Bifurcation  structure in $\alpha$}
We are interested in the behaviour of system \eqref{new} as we vary the value of $\alpha$. For $\alpha$ greater than the first fold of the quartic (the rightmost in Fig.~\ref{bifupper1}), the corresponding equilibrium on the quartic is stable.

\begin{figure}[h!] \centering
  \begin{subfigure}[b]{0.49\linewidth}
    \centering
 \begin{tikzpicture}[thick,scale=0.8, every node/.style={scale=0.8}]
   \node at (0,0){\includegraphics[width=\linewidth]{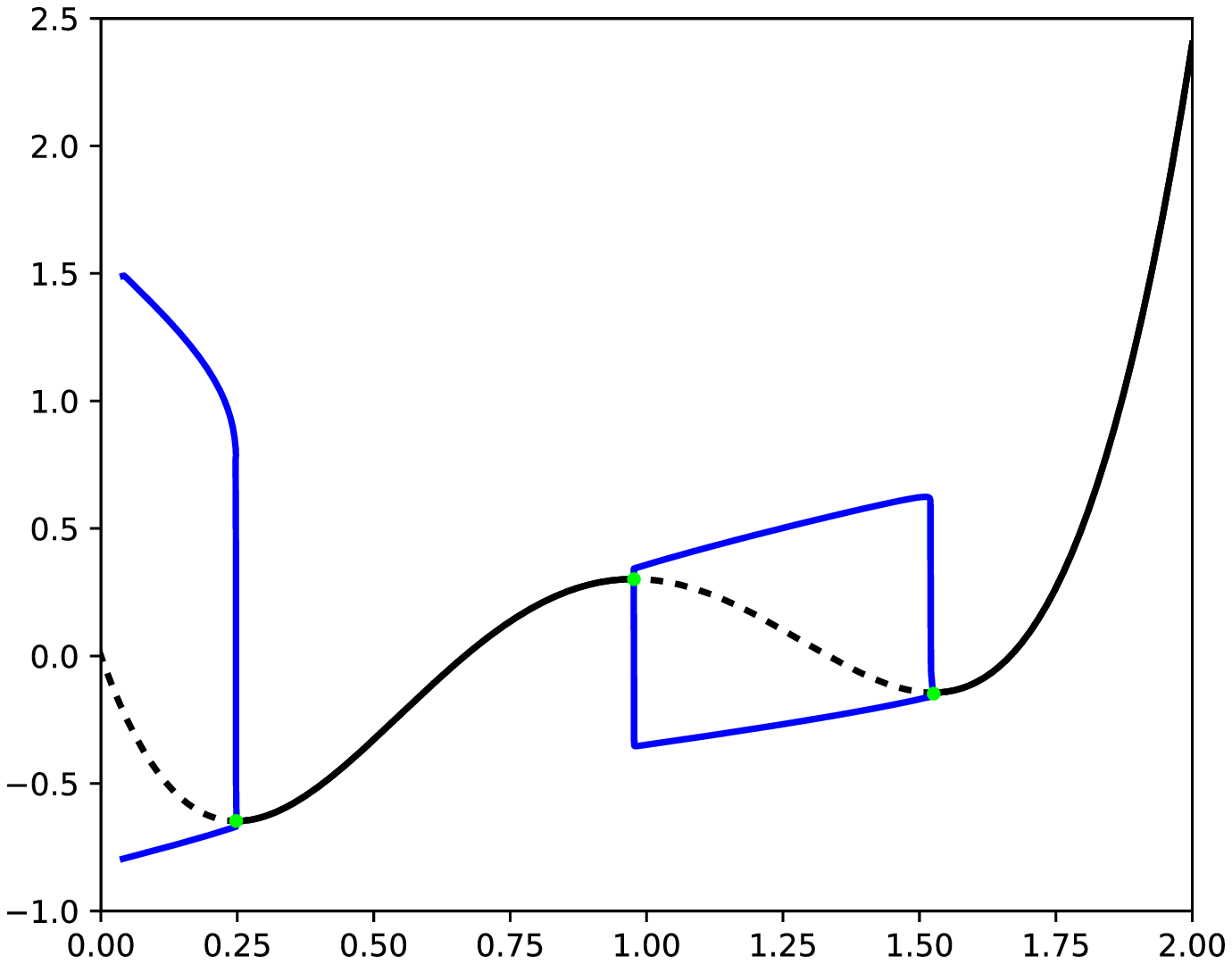}};	
   \node at (0,-3.25) {$\alpha$};
   \node at (-3.5,0.2) {$p_1$};
    \node at (2,-1.7) {\textcolor{green}{$H_1$}};
  \node at (-0.1,-0.4) {\textcolor{green}{$H_2$}};
    \node at (-2,-2.4) {\textcolor{green}{$H_3$}};
	\end{tikzpicture}  
    \caption{} 
    \label{biflower:total} 
  \end{subfigure}
  \begin{subfigure}[b]{0.49\linewidth}
    \centering
\begin{tikzpicture}[thick,scale=0.8, every node/.style={scale=0.8}]
   \node at (0,0){\includegraphics[width=\linewidth]{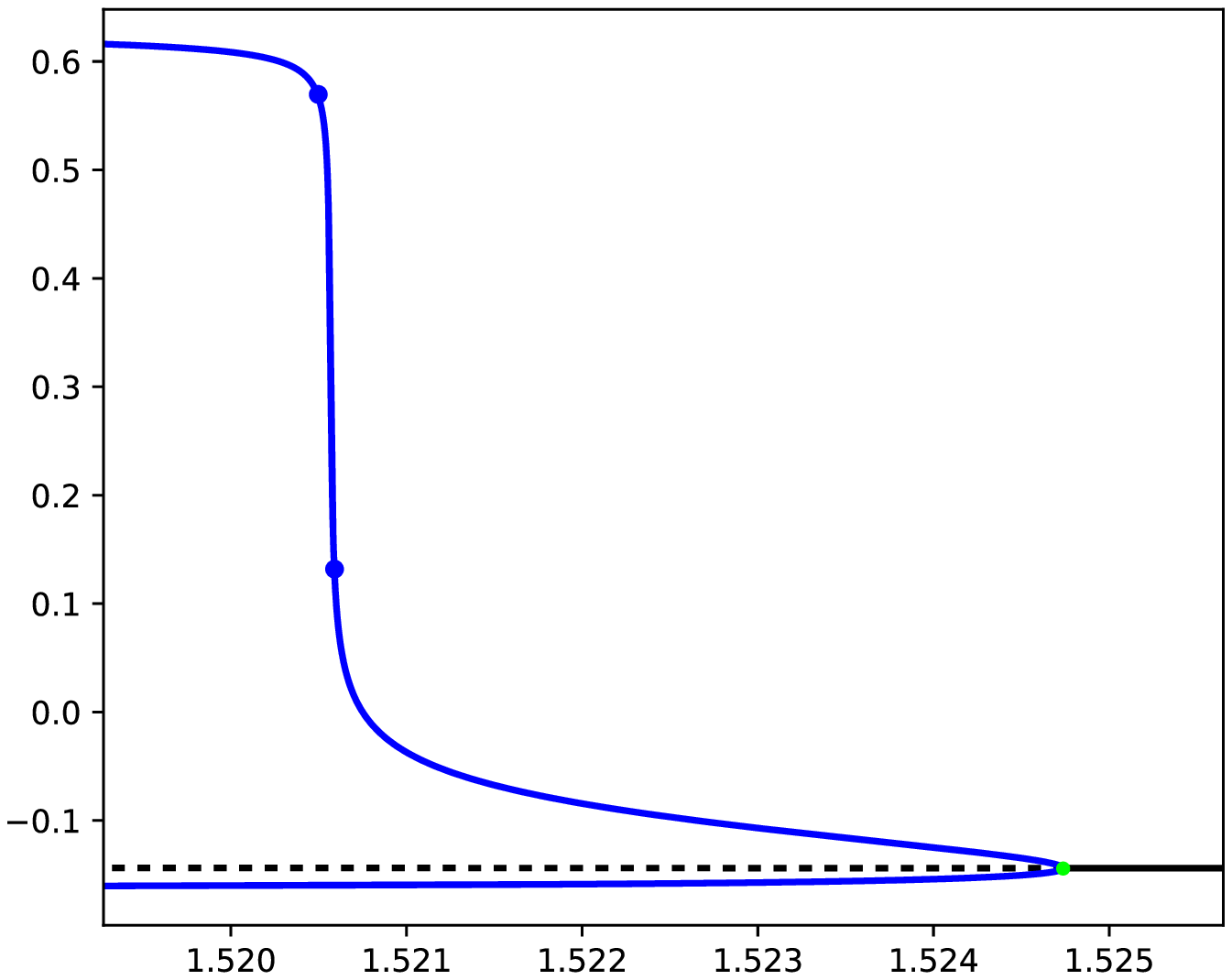}};	
   \node at (0,-3.25) {$\alpha$};
  \node at (-3.5,0.2) {$p_1$};
   \node at (2.5,-2) {\textcolor{green}{$H_1$}};
   \node at (-1.45,-0.5) {\textcolor{blue}{$1$}};
   \node at (-1.5,2.2) {\textcolor{blue}{$2$}};
	\end{tikzpicture}  
    \caption{} 
    \label{biflower:hopf1} 
  \end{subfigure} 
  \begin{subfigure}[b]{0.49\linewidth}
    \centering
\begin{tikzpicture}[thick,scale=0.8, every node/.style={scale=0.8}]
   \node at (0,0){\includegraphics[width=\linewidth]{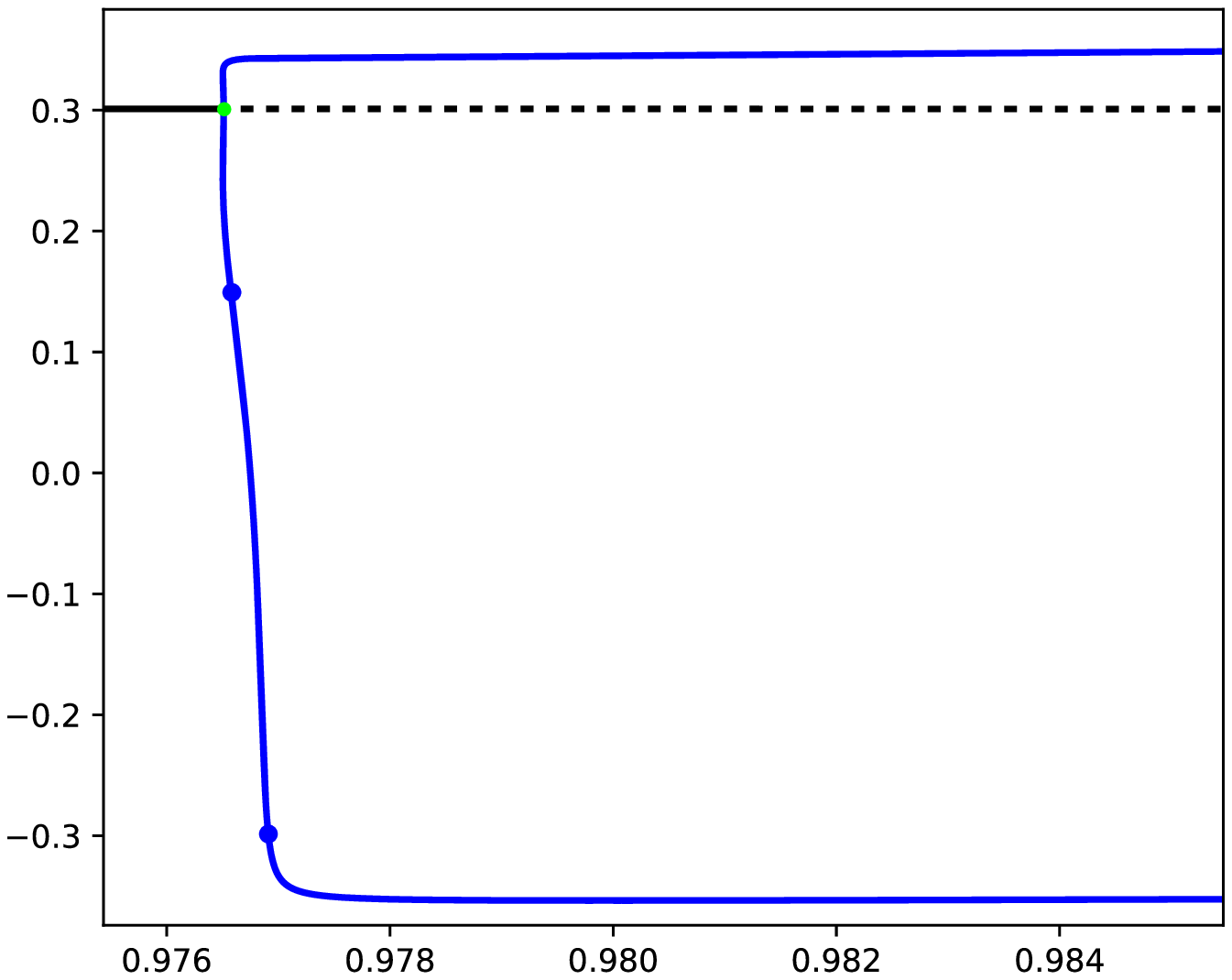}};	
   \node at (0,-3.25) {$\alpha$};
  \node at (-3.5,0.2) {$p_1$};
    \node at (-2,1.8) {\textcolor{green}{$H_2$}};
  \node at (-1.8,-2.1) {\textcolor{blue}{$3$}};
   \node at (-2,1) {\textcolor{blue}{$4$}};
	\end{tikzpicture}  
   \caption{} 
    \label{biflower:hopf2} 
  \end{subfigure}
  \begin{subfigure}[b]{0.49\linewidth}
    \centering
       \begin{tikzpicture}[thick,scale=0.8, every node/.style={scale=0.8}]
   \node at (0,0){\includegraphics[width=\linewidth]{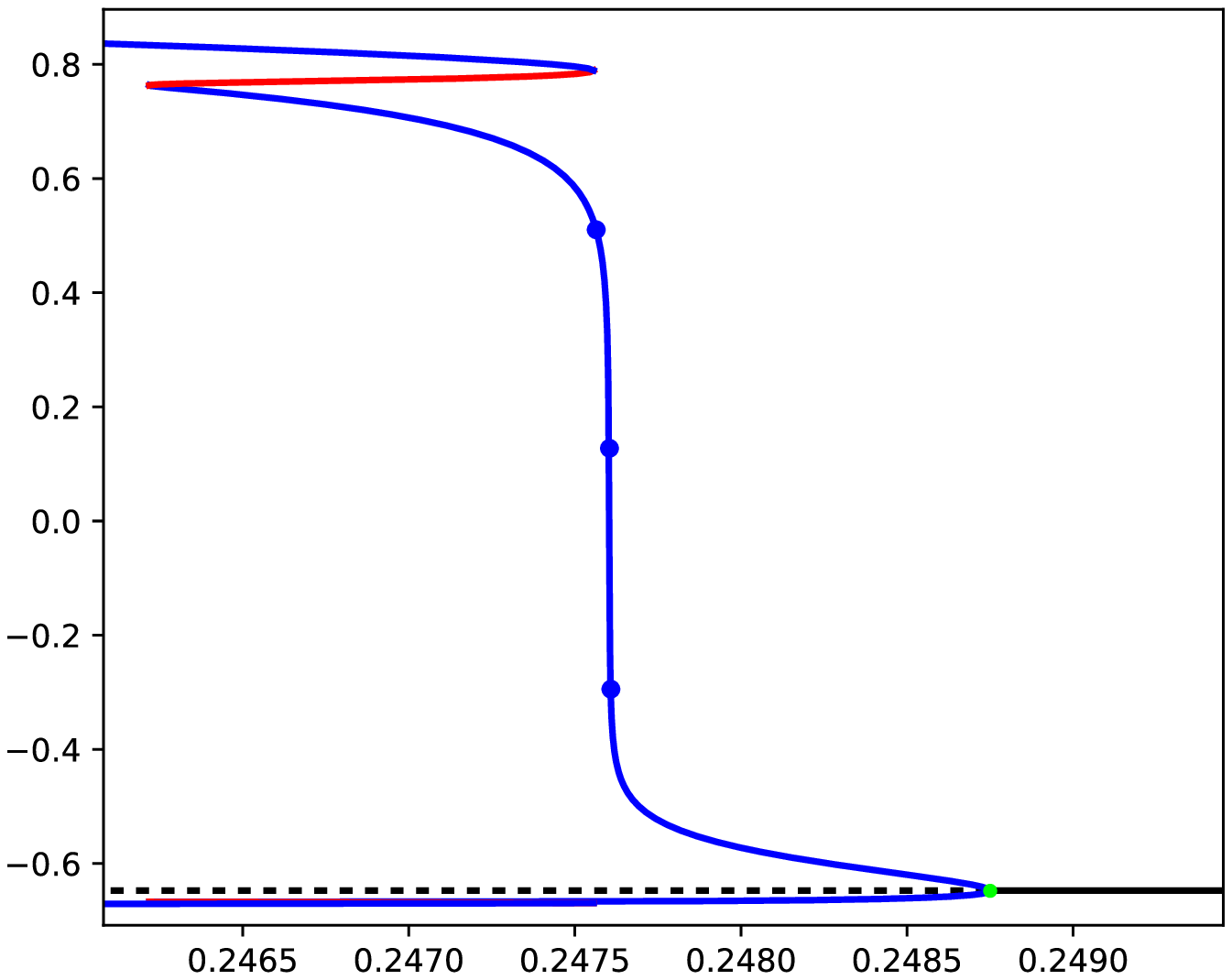}};	
   \node at (0,-3.25) {$\alpha$};
   \node at (-3.5,0.2) {$p_1$};
   \node at (2.1,-2.2) {\textcolor{green}{$H_3$}};
   \node at (0.2,-1.3) {\textcolor{blue}{$5$}};
   \node at (0.15,0.2) {\textcolor{blue}{$6$}};
   \node at (0.1,1.5) {\textcolor{blue}{$7$}};
	\end{tikzpicture}  
    \caption{} 
   \label{biflower:hopf3} 
 \end{subfigure} 
  \caption{(a): Bifurcation diagram of the $(p_1,p_2)$ system with respect to $\alpha$, showing branches of equilibria represented by their $p_1$ value, and branches of limit cycles for which both maximum and minimum of $p_1$ along the cycles, are plotted. Solid black lines correspond to stable equilibria, while dashed black lines correspond to unstable ones. Stable limit cycle branches are shown in blue, and unstable ones in red. (b), (c), (d): zoomed view of the canard-explosive branches emanating from each of the three Hopf bifurcation shown (black dots) in panel (a), from right to left. The numbered dots correspond to canard orbits shown in Fig.~\ref{canardslower}.}
  \label{biflower} 
\end{figure}

At $\alpha=\alpha(H_1)$, we observe a canard explosion; this branch connects to the second Hopf. The branch stemming from the third and last Hopf point on the quartic extends towards $\alpha=0$. When comparing Fig.~\ref{bifupper1} with Fig.~\ref{biflower}, we can observe that the other branch now approaches $\alpha=0$. Hence, for some value of $r_1$ between $5$ and $6.4$ we must observe an exchange of positions within these branches. In Fig.~\ref{bifcon1}, we show a good approximation of the value of $r_1$ at which this exchange happens.

\begin{figure}[h!]\centering
\begin{subfigure}[b]{0.49\linewidth}
   \centering
        \begin{tikzpicture}[thick,scale=0.8, every node/.style={scale=0.8}]
   \node at (0,0){\includegraphics[width=\linewidth]{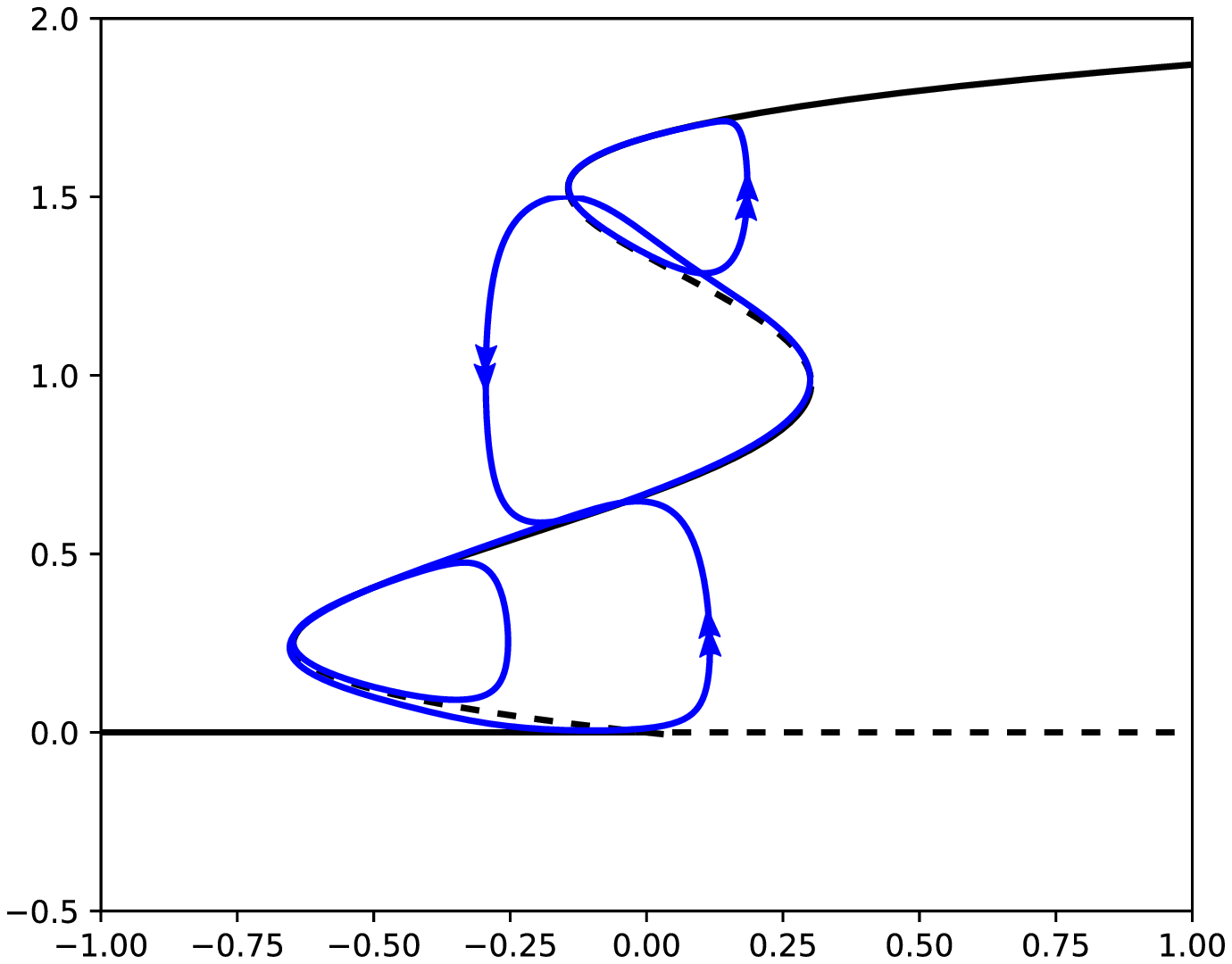}};	
   \node at (0,-3.25) {$p_1$};
   \node at (-3.5,0.25) {$p_2$};
   \node at (1.05,1.5) {\textcolor{blue}{$1$}};
   \node at (-0.9,0.6) {\textcolor{blue}{$3$}};
  \node at (-0.45,-0.75) {\textcolor{blue}{$5$}};
  \node at (0.8,-0.75) {\textcolor{blue}{$6$}};
   \end{tikzpicture}
    \caption{} 
    \label{canlower1} 
  \end{subfigure}
  \begin{subfigure}[b]{0.49\linewidth}
    \centering
\begin{tikzpicture}[thick,scale=0.8, every node/.style={scale=0.8}]
   \node at (0,0){\includegraphics[width=\linewidth]{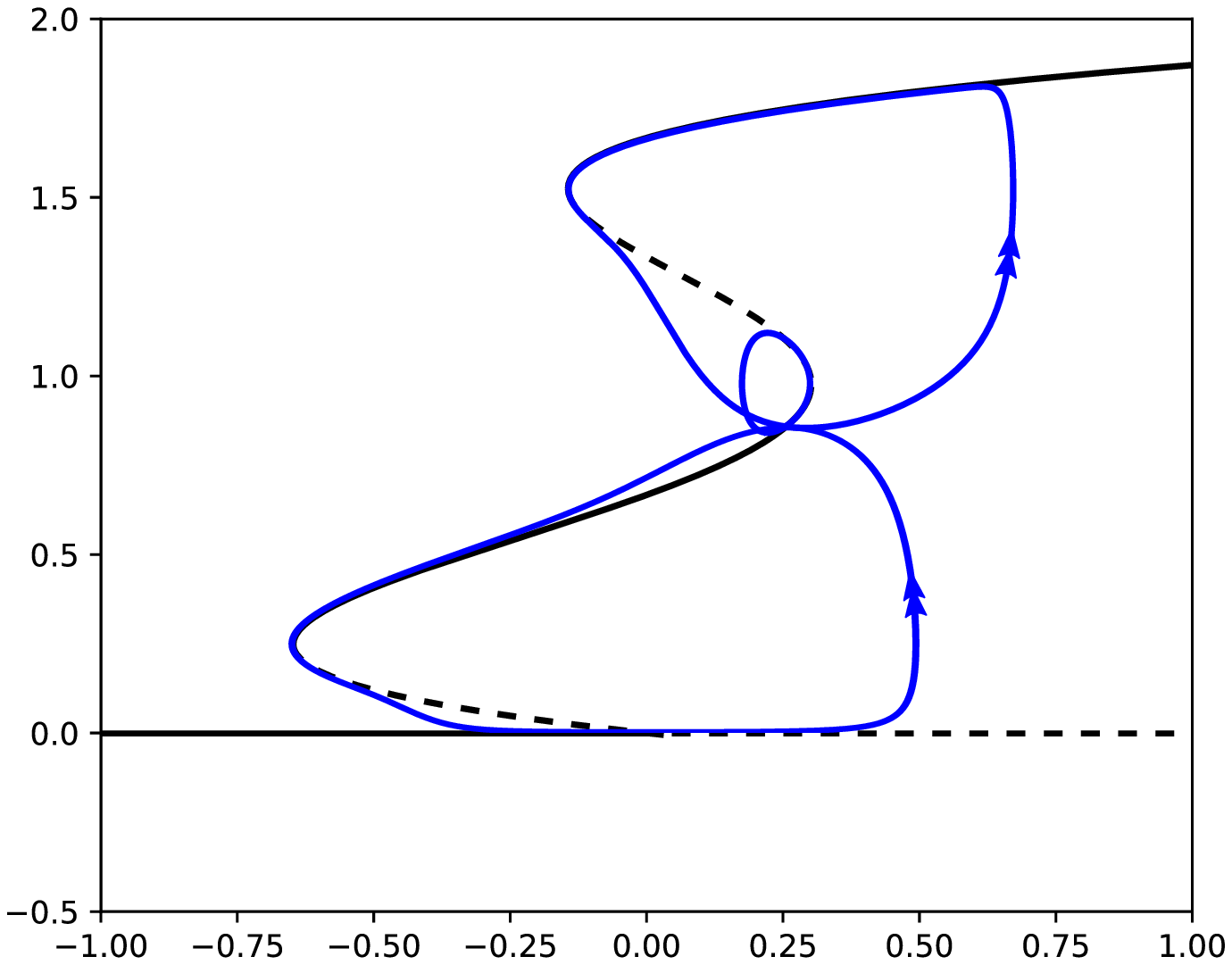}};	
   \node at (0,-3.25) {$p_1$};
   \node at (-3.5,0.25) {$p_2$};
   \node at (2.6,1.5) {\textcolor{blue}{$2$}};
  \node at (1.4,0.6) {\textcolor{blue}{$4$}};
  \node at (2.1,-0.75) {\textcolor{blue}{$7$}};
   \end{tikzpicture}
    \caption{} 
    \label{canlower2} 
 \end{subfigure} 
\caption{Selection of canard cycles that exist along the branches of limit cycles shown in Figs.~\ref{biflower} (b), (c) and (d).}
\label{canardslower}
\end{figure}

In Fig.~\ref{biflower} panels (b-d), we highlight with dots, numbered from 1 to 7, points on the bifurcation diagram corresponding to a selection of canard orbits. These orbits are then plotted in phase plane in Fig.~\ref{canardslower}. These cycles are numbered in decreasing values of $\alpha$, although clearly the cycles belonging to the same canard explosive branch correspond to values of $\alpha$ extremely close.

\subsection{Close to threshold bifurcation diagram}
\begin{figure}[h!]\centering
  \begin{subfigure}[b]{0.49\linewidth}
    \centering
\begin{tikzpicture}[thick,scale=0.8, every node/.style={scale=0.8}]
   \node at (0,0){\includegraphics[width=\linewidth]{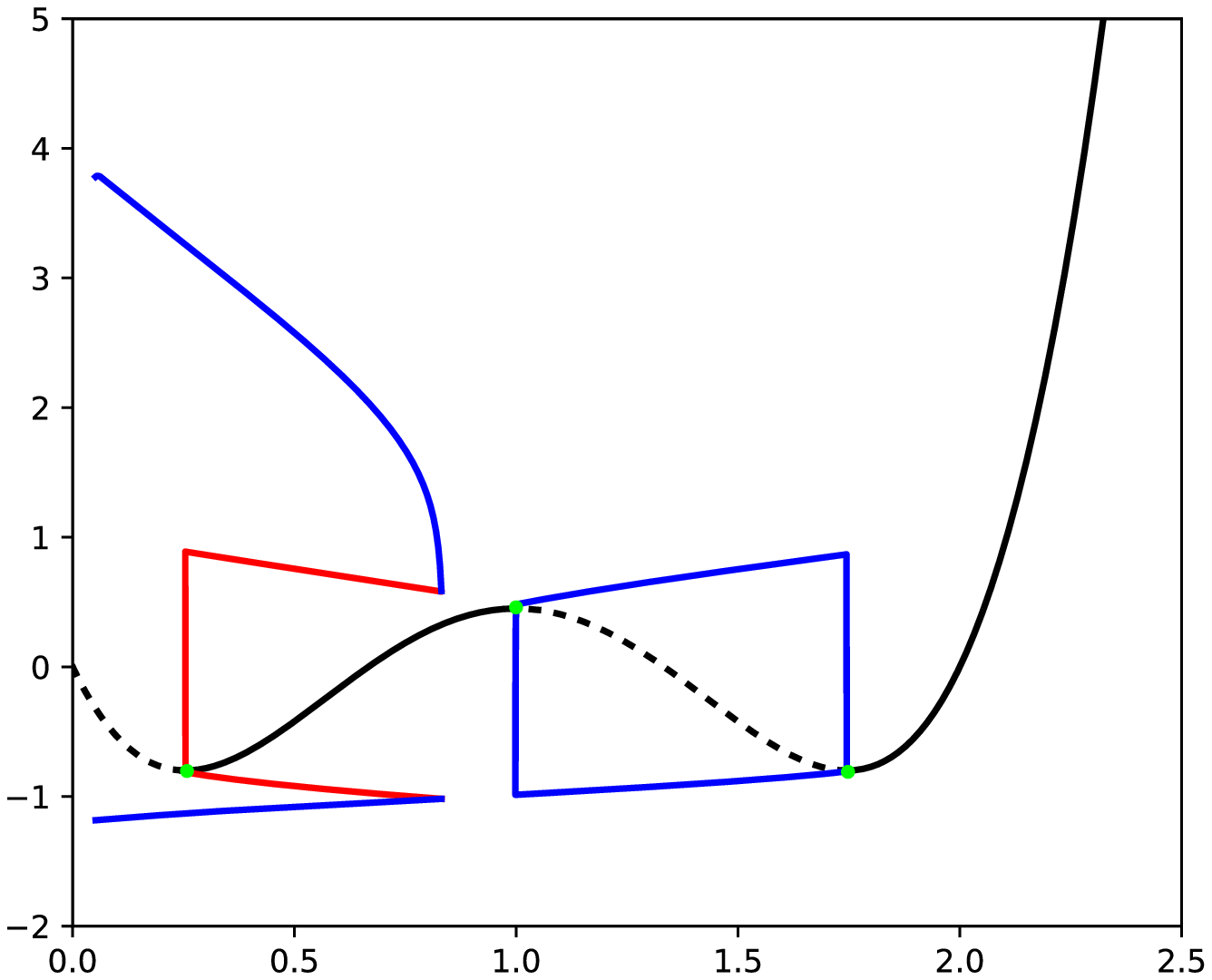}};	
  \node at (0,-3.25) {$\alpha$};
   \node at (-3.5,0.2) {$p_1$};
    \node at (1.8,-2) {\textcolor{green}{$H_1$}};
  \node at (-0.3,-0.5) {\textcolor{green}{$H_2$}};
    \node at (-2.1,-1.4) {\textcolor{green}{$H_3$}};
	\end{tikzpicture}  
    \caption{} 
    \label{bifsym} 
  \end{subfigure}
    \begin{subfigure}[b]{0.49\linewidth}
    \centering
\begin{tikzpicture}[thick,scale=0.8, every node/.style={scale=0.8}]
   \node at (0,0){\includegraphics[width=\linewidth]{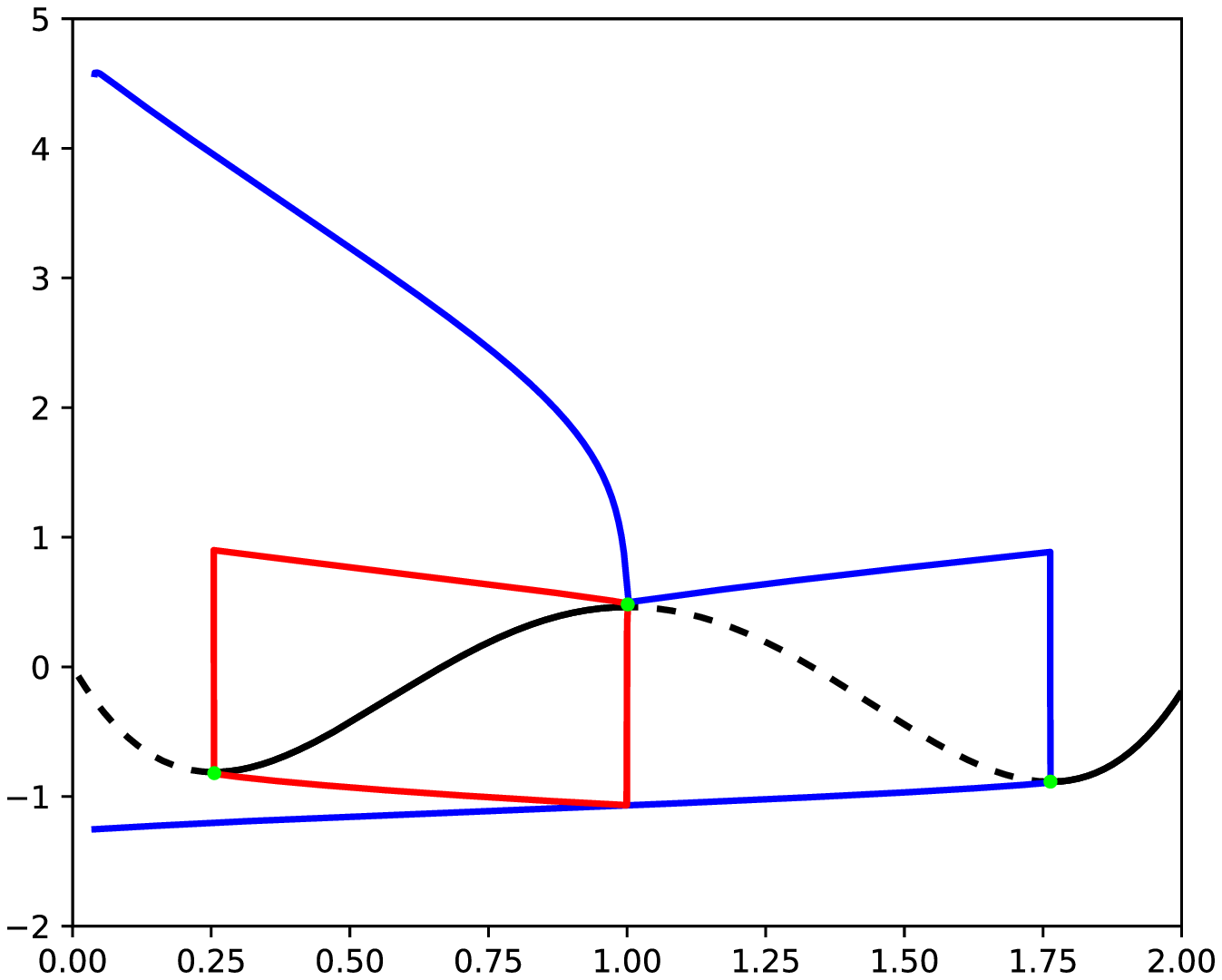}};	
   \node at (0,-3.25) {$\alpha$};
   \node at (-3.5,0.2) {$p_1$};
    \node at (2.6,-2.1) {\textcolor{green}{$H_1$}};
   \node at (0,-0.5) {\textcolor{green}{$H_2$}};
    \node at (-1.9,-1.4) {\textcolor{green}{$H_3$}};
	\end{tikzpicture}  
    \caption{} 
    \label{bif08} 
 \end{subfigure} 
 \begin{subfigure}[b]{0.49\linewidth}
    \centering
\begin{tikzpicture}[thick,scale=0.8, every node/.style={scale=0.8}]
   \node at (0,0){\includegraphics[width=\linewidth]{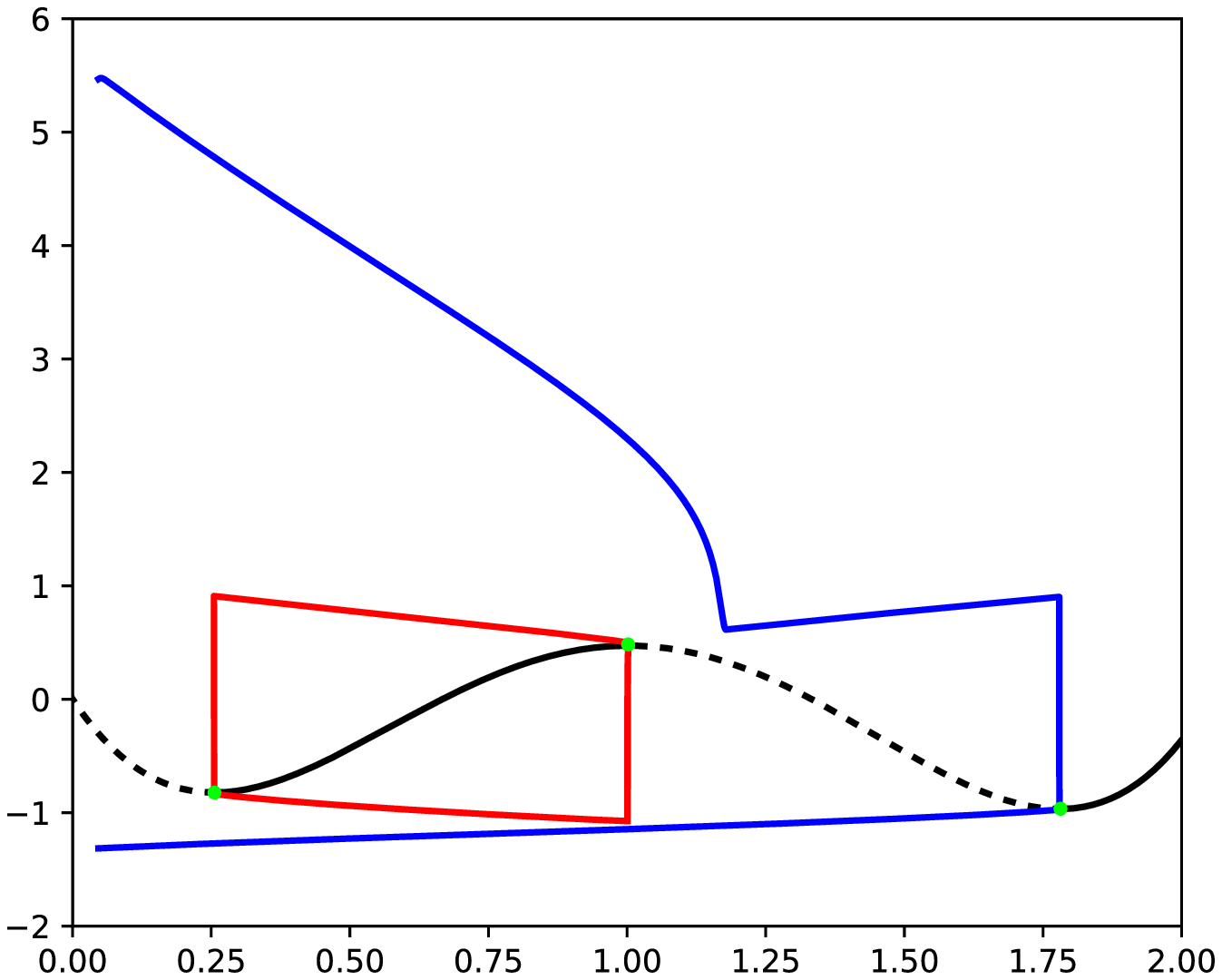}};	
   \node at (0,-3.25) {$\alpha$};
   \node at (-3.5,0.2) {$p_1$};
    \node at (2.8,-2.3) {\textcolor{green}{$H_1$}};
   \node at (0.3,-0.75) {\textcolor{green}{$H_2$}};
    \node at (-1.9,-1.6) {\textcolor{green}{$H_3$}};
	\end{tikzpicture}  
    \caption{} 
    \label{bifother} 
  \end{subfigure} 
\caption{Bifurcation diagram of system~\eqref{new_slow_var} with respect to parameter $\alpha$ and for the three different values of $r_1$ mentioned above, namely: (a) $r_1=6$; (b) $r_1=6.08$; (c) $r_1=6.15$. In all panels, branches of equilibria are shown in their $p_1$ projection and branches of limit cycles are shown in their $(\max(p_1), \min(p_1))$ projection. All branches are explosive and contain canard cycles.}  \label{bifcon1}
\end{figure}

In this section, we showcase the bifurcation diagram of $\alpha$ for system \eqref{new} with the following choices of the quartic equation:
$$
p_1=-(0.05(-3p_2+6\textcolor{white}{.00})(-3p_2+4)(-3p_2+2)3p_2)\quad\text{in Fig.~\ref{bifsym}},
$$
$$
p_1=-(0.05(-3p_2+6.08)(-3p_2+4)(-3p_2+2)3p_2)\quad\text{in Fig.~\ref{bif08}},
$$
$$
p_1=-(0.05(-3p_2+6.15)(-3p_2+4)(-3p_2+2)3p_2)\quad\text{in Fig.~\ref{bifother}}.
$$
When comparing the first and the last of the abovementioned figures, we can conclude that by a continuity argument in $r_1$, there must be a value $r_1^*\in[6, 6.15]$ for which the branches of limit cycles collide, before the connection among the Hopf points changes. Due to the stiffness of the problem for $\eps \ll 1$ and the canard explosions at the three Hopf points, which make numerical exploration even more susceptible to small changes in the parameters, finding the exact value of $r_1$ at which this occurs is quite challenging. Strictly speaking, one could continue the fold of cycles on the branch coming from the leftmost Hopf point in panel (a) and the fold of cycles on the branch coming from the rightmost Hopf point in panel (b), in $(\alpha,r_1)$. However, the stiffness of the problem makes it very hard already to detect reliably these fold bifurcations of cycles -- spurious bifurcations are likely to be flagged by \textsc{auto} --, let alone continuing them in two parameters. If this computation were possible and reliable, then they would collide with the Hopf point at the transcritical crossing that we are expecting for $r_1\in[6,6.15]$.
%

%

\section{Conclusions}
\label{sec:conclusion}
In this article, we have analysed with a mix of analytical and numerical tools a slow-fast dynamical system with a quartic critical manifold, which is pertinent to the modeling of neurotransmitter release. The model is a phenomenological extension of our initial quadratic model from~\cite{rodrigues2016time}, which was calibrated and fitted to experimental data from both excitatory and inhibitory synapses (once we append to it the Tsodyks-Markram synaptic ressources model), and focused on the exocytotic part of the neurotransmitter vesicle cycle. In the present case, we have introduced an extension of this model, which constitutes a theoretical first attempt to incorporate the endocytotic part of the vesicle cycle and capture further multi-timescale effects. We have also explored with more analytical details the generating mechanism behind the delay to the neurotransmitter release via the entry-exit function. However, the difficulty of the data related to endocytosis led us to first propose a theoretical analysis of this extended model. Namely, we have focused on two main aspects: first, the transient response of this quartic model to a time-dependent forcing term akin to an input spike; second, the model's asymptotic behaviour and associated bifurcation structure with respect to a few key parameters.

In terms of response of the model to input spikes, we have identified several interesting configurations of the quartic critical manifold giving rise to interesting responses that were not possible with the quadratic (exocytotic) previous model. In particular, as shown in Fig.~\ref{conveq1}, we can obtain a transient response with a small cycle followed by a large one. As an interesting question for future work, we plan to relate this dynamics of the quartic model to the differential vesicle refiling due to multi-timescale effects during endocytosis, in link with the \textit{kiss-and-run} vs. \textit{kiss-and-stay} scenarios~\cite{watanabe2013,watanabe2015}. To make the present quartic model interpretable in terms of this fine biological phenomenon is a non-trivial task, starting with the meaning of variables, which will have to change from the exocytotic to the endocytotic part of the vesicle cycle.

Regarding the long-term dynamics, this quartic model can exhibit a complex periodic regime with various families of limit cycles that are shaped and organised in parameter space in tight connection with the geometry of the quartic critical manifold. Hence, the second part of the present work was focused on analysing the limit cycle structure of the model depending on the relative position of the fold points of the critical manifold. As explained already in the introduction, the multi-timescale nature of the model, together with the geometry of its nullclines, make it an example of canard system; note that canard cycles were already possible in the quadratic model. However, the quartic geometry allows for multiple families of canard cycles organised in parameter space along explosive branches born at three singular Hopf bifurcations, located near each of the three fold points of the critical manifold. As we have showcased, depending on the relative positions of these three fold points, two limit-cycle branches are connected in parameter space, while the third one terminates at a saddle homoclinic bifurcation. We provided precise numerical evidence of the transition between one possible configuration and the other one.

Throughout the paper, we encountered various interesting analytical questions which are beyond the scope of this paper, but which can inspire future research. In particular, one could investigate how to overcome the limitation we encountered in our bifurcation analysis. Namely, is there a way to complete the numerical continuation of limit cycles converging to a homoclinic or heteroclinic orbit while including a segment that is exponentially close to a flat component of the critical manifold? As the period approaches infinity, the  computational scheme ceases to converge, and one would probably need much more sophisticated numerical methods in order to tackle this challenge.

From a purely analytical point of view, the system of interest is non-autonomous and non-Lipschitz, as a consequence of the input function $V_{\text{in}}$. One could regularise the input function $V_{\text{in}}$ and suspend the resulting smooth equations so as to recover an autonomous system. This would both allow for further analytical treatment of the problem, as well as, enable a generalization to a network setting, in which multiple neurons stimulate each other in a more straightforward way.

Finally, from a GSPT point of view, it would be particularly interesting to provide a general theory of entry-exit functions when the applicable region (i.e. the region in which the slow variable is increasing) disappears as a parameter approaches $0$ ($\alpha$, in our case). In addition to this, it could be interesting to formulate a general theory for a folded critical manifold, such as the quadratic and quartic curves studied in~\cite{rodrigues2016time} and here, respectively, but with more general functions and some base assumption on their shape and behaviour. These are substantial mathematically oriented questions that we leave as possible avenues for follow-ups on this problem.

\bibliographystyle{plain}
\bibliography{biblio}

\end{document}